\documentclass[10pt]{amsart}
\usepackage[english]{babel}
\usepackage{graphicx}
\usepackage{amsmath,amssymb,hyperref}
\newcommand{\diff}[2]{\mbox{{\rm Diff}{${\,}_{#1}({\mathbb C}^{#2},0)$}}}
\newcommand{\diffh}[2]{\mbox{$\widehat{\rm Diff}{{\,}_{#1}({\mathbb C}^{#2},0)}$}}

\newcommand{\cn}[1]{\mbox{(${\mathbb C}^{#1},0$)}}

\newcommand{\ex}{\'{e}}

\newcommand{\ox}{\'{o}}

\newtheorem{pro}{Proposition}[section]
\newtheorem{teo}{Theorem}[section]

\newtheorem{cor}{Corollary}[section]

\newtheorem{lem}{Lemma}[section]
\newtheorem{rem}{Remark}[section]
\newtheorem{defi}{Definition}[section]
\begin{document}

\title[Formal classification of up-diffeomorphisms]
{Formal classification of  unipotent parameterized diffeomorphisms}


\author{Javier Rib\ox n}
\thanks{Universidad de Valladolid,
Departamento de Algebra, Geometria y Topologia, Paseo de la Magdalena s/n, Valladolid, Spain, 47011}
\thanks{e-mail address: jfribon@impa.br}
\thanks{MSC-class: 37F75, 32H02, 32S65, 32A05}
\date{\today}
\maketitle

\bibliographystyle{plain}
\section*{Abstract}
  We provide a complete system of invariants for the formal classification
of complex analytic unipotent germs of diffeomorphism at $\cn{n}$ fixing
the orbits of a regular vector field.
%
%
We reduce the formal classification problem to solve a linear differential
equation. Then we describe the formal invariants; their nature depends on the
position of the fixed points set $Fix \varphi$ with respect to the regular vector field preserved
by $\varphi$. We get invariants specifically attached to higher dimension ($n \geq 3$)
although generically they are analogous to the one-dimensional ones.
\section{Introduction}
  We provide a complete system of invariants for the formal classification
of complex analytic unipotent parameterized germs of diffeomorphism at $\cn{n+1}$.
Consider coordinates $(x,x_{1},\hdots,x_{n})$ in ${\mathbb C}^{n+1}$.
Denote by $\diff{}{n+1}$ the group of complex analytic germs of diffeomorphism at $\cn{n+1}$.
We define the group
\[ \diff{p}{n+1} = \{ \varphi \in \diff{}{n+1} : x_{j} \circ \varphi = x_{j} \ {\rm for \ all} \ 1 \leq j \leq n\}  \]
of parameterized diffeomorphisms. The group
\[ \diff{up}{n+1} =
\left\{{ \varphi \in \diff{p}{n+1} : \frac{\partial{(x \circ \varphi)}}{\partial{x}}
(0,\hdots,0) = 1}\right\}   \]
is the set of unipotent elements of $\diff{p}{n+1}$. If $\varphi \in \diff{up}{n+1}$
then we say that $\varphi$ is a unipotent parameterized diffeomorphism
(or up-diffeomorphism for shortness). Let $\varphi \in \diff{p}{n+1}$, we have
that $(\partial{(x \circ \varphi)}/\partial{x})(0)=1$ if and only if
$x \circ \varphi(x,0,\hdots,0) \in \diff{}{}$ is tangent to the identity.
Thus the set of up-diffeomorphisms is the set of finite dimensional holomorphic
perturbations of tangent to the identity germs of diffeomorphism.

  The complex analytic germs of diffeomorphism in one complex variable
are well-known. Those germs whose linear part is not periodic are formally
linearizable. On the one hand they are analytically linearizable if the linear part
is not a rotation. On the other hand we find ``small divisor problems"
(Siegel \cite{Siegel}, Bruno \cite{Bru}, Yoccoz \cite{Yo},
P\ex rez-Marco \cite{PM:acta}) leading to very complicated dynamics if
we deal with non-linearizable diffeomorphisms whose fixed point is of indifferent type.

 The study of the diffeomorphisms with periodic linear part
can be reduced to the one of tangent to the identity diffeomorphisms
where we know the formal, topological (Leau \cite{Leau}, Camacho \cite{Cam}) and
analytical (Birkhoff, Ecalle, Voronin \cite{V}, Malgrange \cite{mal:ast}) classifications.
The only topological invariant is the order of contact with the identity;
this discrete invariant plus a numerical invariant called residue
(cf. subsection \ref{subsec:narefu}) compose a complete system of formal invariants.
The analytical classification is more complicated; we can express
the invariants as a collection (changes of chart) of one-variable germs of
diffeomorphism (Martinet and Ramis \cite{MaRa:aen}). The number of changes of chart is twice
the order of contact with the identity.

A natural generalization of germs of diffeomorphism at $\cn{}$ are parameterized
diffeomorphisms. We are interested on the formal classification of parameterized diffeomorphisms.
We denote the fixed points set of a diffeomorphism $\varphi$ by $Fix \varphi$.
Consider $\varphi \in \diff{p}{n+1}$ such that
$(\partial{(x \circ \varphi)}/\partial{x})(0)$ is not a root of the unit.
The function linear part
$\partial{(x \circ \varphi)}/\partial{x}:Fix \varphi \to {\mathbb C}$
is the only formal invariant attached to $\varphi$ as in the one-dimensional case.
Thus the task of
obtaining a formal classification in $\diff{p}{n+1}$ can be reduced
to exhibit a complete system of formal invariants for up-diffeomorphisms.

  We denote by $\diffh{}{n+1}$, $\diffh{p}{n+1}$ and $\diffh{up}{n+1}$ the formal completions
of $\diff{}{n+1}$, $\diff{p}{n+1}$ and $\diff{up}{n+1}$ respectively.

 A unipotent $\varphi \in \diff{}{n+1}$ is the exponential of a unique formal nilpotent vector field
(see section \ref{sec:baproup} for definitions) that we denote by $\log \varphi$.
Consider $\varphi$ in $\diff{up}{n+1}$;
we have that $\log \varphi$ is of the form ${\rm exp}(\hat{u} (x \circ \varphi- x) \partial/\partial{x})$
where $\hat{u} \in {\mathbb C}[[x,x_{1},\hdots,x_{n}]]$ is a unit. The logarithm of $\varphi$
can be extended to $Fix \varphi$, more precisely
\begin{pro}
\label{pro:iforco}
Let $\varphi = {\rm exp}(\hat{u} (x \circ \varphi- x) \partial/\partial{x}) \in \diff{up}{n+1}$. Then
$\hat{u}$ belongs to $\lim_{\leftarrow} {\mathbb C}\{x,x_{1},\hdots,x_{n}\}/I(x \circ \varphi -x)^{j}$.
\end{pro}
In other words there exists $u_{j} \in {\mathbb C}\{x,x_{1},\hdots,x_{n}\}$
such that $\hat{u} -u_{j} \in (x \circ \varphi -x)^{j}$ for all $j \in {\mathbb N}$.

We say that a germ of analytic variety at $({\mathbb C}^{n+1},0)$ is {\it fibered} if it is
a union of orbits of $\partial{}/\partial{x}$. By definition
$\hat{\rho} \in \diffh{}{n+1}$ is special with respect to $f=0$ (for $f \in {\mathbb C}\{x,x_{1},\hdots,x_{n}\}$)
if $\hat{\rho} \in \diffh{p}{n+1}$ and
$\hat{\rho}_{|\gamma} \equiv Id \Leftrightarrow x \circ \hat{\rho} - x \in I(\gamma)$ for all
non-fibered irreducible component $\gamma$ of $f=0$.
\begin{pro}
Let $\varphi_{1},\varphi_{2} \in \diff{up}{n+1}$. Assume that $\varphi_{1}$ and $\varphi_{2}$
are formally conjugated. Then there exists $\sigma \in \diff{}{n+1}$ and a special
$\hat{\sigma} \in \diffh{p}{n+1}$ (with respect to $x \circ \varphi_{2} - x=0$) such that
$(\hat{\sigma} \circ \sigma) \circ \varphi_{1} = \varphi_{2} \circ (\hat{\sigma} \circ \sigma)$.
\end{pro}
The last proposition implies that up to analytic change of coordinates
every couple of formally conjugated up-diffeomorphisms
are conjugated by a special element of $\diffh{p}{n+1}$. We study
the equivalence relation in $\diff{up}{n+1}$ given by
$\varphi_{1} \sim \varphi_{2}$ if $\varphi_{1}$ and $\varphi_{2}$ are formally
conjugated by a special transformation with respect to $Fix \varphi_{1}$.
Every class of equivalence is contained in a set
\[ {\mathcal D}_{f} = \{ \varphi \in \diff{up}{n+1} : (x \circ \varphi - x)/f \ {\rm is \ a \ unit} \} . \]
The classes of this equivalence relation are connected sets in the compact-open topology.
As a consequence to determine whether or not there exists a formal special conjugation between
up-diffeomorphisms can be reduced to solve a linear problem. More precisely we can associate to
$\varphi_{1},   \varphi_{2} \in {\mathcal D}_{f}$ the homological equation
\[ \frac{\partial{\alpha}}{\partial{x}} = \frac{1}{f}
\left({
\frac{1}{\hat{u}_{1}} - \frac{1}{\hat{u}_{2}}
}\right)  \]
where $\hat{u}_{j}=(\log \varphi_{j})(x)/f$ for $j \in \{1,2\}$.
Let $\prod_{j=1}^{p} f_{j}^{l_{j}} \prod_{j=1}^{1} F_{j}^{m_{j}}$ be the decomposition of $f$
in irreducible components. By choice $F_{j}=0$ is fibered for $1 \leq j \leq q$ whereas
$f_{k}=0$ is non-fibered for $1 \leq k \leq p$. We say that the homological equation is special
(with respect to $f$) if there exists a solution of the form
$\alpha = \hat{\beta}/(\prod_{j=1}^{p} f_{j}^{l_{j}-1} \prod_{j=1}^{q} F_{j}^{m_{j}})$
where $\hat{\beta} \in {\mathbb C}[[x,x_{1},\hdots,x_{n}]]$. Such a solution is also called special.
We have
\begin{pro}
Let $\varphi_{1},\varphi_{2} \in {\mathcal D}_{f} \subset \diff{up}{n+1}$. Then $\varphi_{1}$ and $\varphi_{2}$
are formally conjugated by a special transformation if and only if
the homological equation associated to $\varphi_{1}$ and $\varphi_{2}$ is special.
\end{pro}
Let $\varphi={\rm exp}(\hat{u} f \partial/\partial{x}) \in {\mathcal D}_{f}$.
The formal 1-form $dx/(\hat{u}f)$ is the dual of $\log \varphi$ in the
relative cohomology of the vector field $\partial / \partial{x}$.
Moreover there exists $u$ in ${\mathbb C}\{x,x_{1},\hdots,x_{n}\}$ such that
$\hat{u}-u \in (f)$ by proposition \ref{pro:iforco}. Therefore we obtain
\[ \frac{dx}{\hat{u}f} - \frac{dx}{uf} = \frac{1}{\hat{u}u} \frac{u-\hat{u}}{f} dx. \]
Since the right hand side does not have poles then the formal properties of
$dx/(\hat{u}f)$ and $dx/(uf)$ are the same. The only formal invariant
of $\varphi \in {\mathcal D}_{f} \subset \diff{up}{}$
for the special conjugation is the
residue of $dx/(uf)$ at $0$. The generalization of this invariant in the higher dimensional case
is the collection of residues of $dx/(uf)$ at $Fix \varphi$.
This collection defines a meromorphic function in every non-fibered irreducible component of $f=0$.

  There are other invariants which are purely related to higher dimension.
For a non-zero $f \in {\mathbb C}\{x,x_{1},\hdots,x_{n}\}$ we define the additive group $Fr(f)$
of homological equations $\partial{\alpha}/\partial{x} = A/f$ ($A \in {\mathbb C}\{x,x_{1},\hdots,x_{n}\}$)
such that $(A/f)dx$ has vanishing residues. Moreover we denote by $Sp(f)$ the subgroup of $Fr(f)$ of
special equations.
\begin{teo}
A complete system of formal invariants for the special conjugation in
${\mathcal D}_{f} \subset \diff{up}{n+1}$ is composed by the residue functions
plus the complex vector space $Fr(f)/Sp(f)$.
\end{teo}
For ${\mathcal D}_{f} \subset \diff{up}{n+1}$ with $n \leq 1$ the only invariants are the residues,
in other words we have $Fr(f)/Sp(f)=0$. The situation is different in higher dimension;
for instance for $f_{0}={(x_{2}-xx_{1})}^{2}$ and $D_{f_{0}} \subset \diff{up}{3}$ we have that
$\dim_{\mathbb C} Fr(f_{0})/Sp(f_{0}) =1$.
Moreover we have $\dim_{\mathbb C} Fr(f)/Sp(f) < + \infty$
for  ${\mathcal D}_{f} \subset \diff{up}{3}$. Thus besides the residue
functions there are only finitely many linear invariants.
In spite of that $Fr(f_{0})/Sp(f_{0}) \sim {\mathbb C}\{x_{3},\hdots,x_{n}\}$
is infinite dimensional
if ${\mathcal D}_{f_{0}}$ is considered as a subset of $\diff{up}{n+1}$ for $n+1 \geq 4$.

 The nature of $Fr(f)/Sp(f)$ depends on the {\it evil set} $S(f)$ of $f$.
This set is the union of the orbits of $\partial/\partial{x}$
contained in non-fibered irreducible components $\gamma$ of $f=0$ such that $f \in I(\gamma)^{2}$.
\begin{pro}
\label{pro:nomoin}
Let $0 \neq f \in {\mathbb C}\{x,x_{1},\hdots,x_{n}\}$ such that $S(f)=\emptyset$. Then
$Fr(f)/Sp(f)=0$.
\end{pro}
Consider the set
\[ K(n)= \{ f \in {\mathbb C}\{x,x_{1},\hdots,x_{n}\} : f(0)=(\partial{f}/\partial{x})(0)=0 \} \]
endowed with the Krull topology. The set
\[ E(n) = \{ f \in K(n) : f(x,0,\hdots,0) \not \equiv 0 \} \]
is open and dense in $K(n)$. Moreover $S(f)=\emptyset$ for all $f \in E(n)$.
\begin{teo}
Fix $n \in {\mathbb N}$.
There exists a dense open subset $E$ of $K(n)$ such that
for $f \in E$ the residue functions provide a complete system of
formal invariants for the special conjugation in ${\mathcal D}_{f} \subset \diff{up}{n+1}$.
\end{teo}
\section{Notations and definitions}
We deal with complex analytic
germs of diffeomorphism defined at $({\mathbb C}^{n+1},0)$. Consider
coordinates $(x,x_{1},\hdots,x_{n})$. We define the group
\[ \diff{p}{n+1} = \{ \varphi \in \diff{}{n+1}  :  x_{j} \circ \varphi = x_{j} \
\forall 1 \leq j \leq n \}  \]
of parameterized diffeomorphisms.
We denote by $\diff{u}{n+1}$ the subgroup of $\diff{}{n+1}$ of unipotent germs of diffeomorphism.
An element $\varphi$ of $\diff{}{n+1}$ is unipotent if its linear part is unipotent, in other
words if $j^{1} \varphi$ has the unique eigenvalue $1$. As a consequence $\varphi \in \diff{p}{n+1}$
is unipotent if and only if $(\partial{(x \circ \varphi)}/\partial{x})(0)=1$.
We will study the elements in the group
\[ \diff{up}{n+1} \stackrel{def}{=} \diff{p}{n+1} \cap \diff{u}{n+1}  \]
 of unipotent parameterized diffeomorphisms.
For the sake of simplicity we will usually replace the expression
unipotent parameterized diffeomorphism with the shorter up-diffeomorphism. The groups that we just
defined have formal completions, we will denote them $\diffh{p}{n+1}$, $\diffh{u}{n+1}$ and
$\diffh{up}{n+1}$.

  The unipotent germs of diffeomorphism are related with nilpotent vector fields. We denote by
${\mathcal X} \cn{n+1}$ the set of germs of complex analytic
vector field which are singular at $0$. We denote
by ${\mathcal X}_{N} \cn{n+1}$ the subset of ${\mathcal X} \cn{n+1}$ of nilpotent vector fields.
The formal completions of these spaces are denoted by
$\hat{\mathcal X} \cn{n+1}$ and $\hat{\mathcal X}_{N} \cn{n+1}$ respectively.
\section{Basic properties of the unipotent parameterized diffeomorphisms}
\label{sec:baproup}
We denote by ${\rm exp}(tX)$ the flow of the vector field $X$, it is the unique solution
of the differential equation
\[ \frac{\partial}{\partial{t}} {\rm exp}(tX) = X({\rm exp}(tX)) \]
with initial condition ${\rm exp}(0X)=Id$. The flow can be developed in power series, such
a property allows to define the formal flow for formal vector fields.
We define $X^{0}(g)=g$ and $X^{j+1}(g)=X^{j}(X(g))$ for all $j \geq 0$. The exponential application
can be expressed in the form
\[ {\rm exp}(tX) = \left({
\sum_{j=0}^{\infty} {t}^{j} \frac{X^{j}(x)}{j!} , \sum_{j=0}^{\infty} {t}^{j} \frac{X^{j}(x_{1})}{j!} ,
\hdots , \sum_{j=0}^{\infty} {t}^{j} \frac{X^{j}(x_{n})}{j!}
}\right) . \]
For any formal nilpotent vector field $\hat{X}$ the sums defining ${\rm exp}(t \hat{X})$ converge
in the Krull topology. The exponential application is by definition ${\rm exp}(1 \hat{X})$.
Next proposition is classical; it relates formal nilpotent vector fields and formal unipotent transformations.
\begin{pro}
The exponential application induces a bijective mapping from
$\hat{\mathcal X}_{N} \cn{n}$ onto $\diffh{u}{n}$. Moreover, if $\hat{X} \in \hat{\mathcal X}_{N} \cn{n}$
then every component of ${\rm exp}(t \hat{X})$ belongs to
${\mathbb C}[t][[x_{1} , \hdots , x_{n}]]$.
\end{pro}
We will call logarithm of a formal up-transformation $\varphi$ the only formal nilpotent vector
field whose exponential is $\varphi$. We denote by $\log \varphi$ the logarithm of $\varphi$.
\subsection{Logarithm of a up-diffeomorphism and the fixed points set}
A up-diffeomorphism preserves the fibration $dx_{1}= \hdots = dx_{n} =0$. Somehow its logarithm
has to preserve the same fibration too.
\begin{pro}
\label{pro:strlog}
  Let $\varphi \in \diffh{u}{n+1}$. Then $\varphi \in \diffh{p}{n+1}$
if and only if $\log \varphi$ can be expressed in the form $\hat{f} \partial / \partial{x}$.
\end{pro}
\begin{rem}
Since the logarithm $\log \varphi = \hat{f} \partial / \partial{x}$ of $\varphi \in \diffh{up}{n+1}$ is
nilpotent then $\hat{f}(0)=0$ and $(\partial{\hat{f}} / \partial{x}) (0) =0$.
\end{rem}
\begin{proof}[proof of proposition \ref{pro:strlog}]
  The implication $(\Leftarrow)$ is a direct consequence of the formula defining ${\rm exp}(\log \varphi)$.

Let us prove the implication $(\Rightarrow)$. The components of
\[ {\rm exp}(t \log \varphi) = ({\varphi}_{0}(t,x,x_{1},\hdots,x_{n}) , \hdots ,
{\varphi}_{n}(t,x,x_{1},\hdots,x_{n})) \]
belong to  ${\mathbb C}[t][[x, x_{1} , \hdots , x_{n}]]$ for $0 \leq j \leq n$.
Since $\diffh{p}{n+1}$ is a group then
${\varphi}_{j}(t,x,x_{1}, \hdots , x_{n}) = x_{j}$
for all $t \in {\mathbb Z}$ and $1 \leq j \leq n$. For $1 \leq j \leq n$ the power series
${\varphi}_{j} - x_{j}$ is identically $0$ because it is $0$ for $t \in {\mathbb Z}$ and it is
polynomial in $t$. Since
\[ (\log \varphi) (x_{j}) = \lim_{t \to 0} \frac{x_{j} \circ {\rm exp}(t \log \varphi) - x_{j}}{t} =
\lim_{t \to 0} 0 = 0 \]
for $1 \leq j \leq n$ then $\log \varphi$ is of the form $\hat{f} \partial{}/\partial{x}$.
\end{proof}
  A up-diffeomorphism $\varphi$ has a fixed points set $x \circ \varphi - x =0$; it is a germ of
hypersurface if $\varphi \neq Id$. We prove next that the fixed points sets of $\varphi$ and
the singular set of $\log \varphi$ coincide.
\begin{pro}
  Let $\varphi={\rm exp}(\hat{f} \partial{}/\partial{x})$ be a up-diffeomorphism. There exists
a formal unit $\hat{u} \in {\mathbb C} [[ x,x_{1},\hdots,x_{n} ]]$ such that
$x \circ \varphi -x = \hat{u} \hat{f}$.
\end{pro}
\begin{proof}
  We define
\[ h_{0} = x , \ h_{1}=\hat{f}, \ \hdots, \ h_{j+1}=\hat{f} \partial{h_{j}}/\partial{x} \ \ \forall j>0. \]
We have $\varphi = (\sum_{j=0}^{\infty} h_{j}/j!,x_{1}, \hdots , x_{n})$.
Since $\hat{f} \partial{}/\partial{x}$ is nilpotent then $h_{j} / \hat{f}$
belongs to the maximal ideal of ${\mathbb C} [[ x,x_{1},\hdots,x_{n} ]]$ for $j \geq 2$. We define
$\hat{u} = 1 + \sum_{j=2}^{\infty} h_{j}/(j! \hat{f})$.
Clearly the series $\hat{u}$ is the unit we were looking for.
\end{proof}
\section{Formal transversality of the logarithm}
  The logarithm of a up-diffeomorphism can be extended to the fixed points set.
Roughly speaking, the logarithm is convergent in the tangent directions to
the fixed points set but it can diverges in the transversal direction. We introduce some
definitions to make rigorous this very simple idea.

  The formal completion of a complex space $(U,\Theta(U))$ ($U$ is a topological space
and $\Theta(U)$ is its sheaf of analytic functions) along a sub-variety $V$ given
by a sheaf of ideals $I$ is the space $(U,\hat{\Theta}_{I}(U))$ where
\[ \hat{\Theta}_{I}(U) = \lim_{\leftarrow} \frac{\Theta(U)}{{I}^{j}} . \]
Throughout this paper we consider three types of formal transversality:
\begin{defi}
  Let $V \subset ({\mathbb C}^{n+1},0)$ be a germ of analytic variety given by an ideal $I(V)$. A series
$\hat{g} \in {\mathbb C}[[x,x_{1},\hdots,x_{n}]]$ is
\begin{itemize}
\item transversally formal along $V$ (or in a equivalent way t.f. along $V$) if
$\hat{g} \in \lim_{\leftarrow} {\mathbb C} \{x,x_{1},\hdots,x_{n} \}/I(V)^{j}$. \\
\item uniformly transversally formal along $V$ (or u.t.f. along $V$) if
$\hat{g}$ belongs to $\lim_{\leftarrow} \vartheta(U)/I(V)^{j}$
for some neighborhood of the origin $U$. \\
\item uniformly semi-meromorphic along $V$ (or u.s.m. along $V$) if
$\hat{g}$ belongs to $\lim_{\leftarrow} {(\vartheta(U))}_{I(V)}/I(V)^{j}$
for some neighborhood of the origin
$U$. Note that ${(\vartheta(U))}_{I(V)}$ is the localized of the
ring of holomorphic functions in $U$ with respect to the ideal $I(V)$.
\end{itemize}
\end{defi}
For the u.s.m. definition we suppose that $V$ is irreducible. In general we say that $\hat{g}$ is
u.s.m. along $V$ if $\hat{g}$ is u.s.m. along every irreducible component of $V$.

The analytic spaces that we complete are
$(0,{\mathbb C} \{x,x_{1},\hdots,x_{n} \})$, $(U,\vartheta(U))$ and $(U,{(\vartheta(U))}_{I(V)})$
respectively. We say that ${X} \in \hat{\mathcal X} \cn{n+1}$ is t.f. along $V$ if
all the functions ${X}(x)$, ${X}(x_{1})$, $\hdots$, ${X}(x_{n})$ are t.f. along $V$.
The other kinds of formal transversality for formal vector fields are defined in an analogous way.

Denote the fixed points set of a diffeomorphism $\varphi$ by $Fix \varphi$.
  Let $\varphi$ be a up-diffeomorphism and consider an irreducible component $\gamma$ of $Fix \varphi$.
We say that $\gamma$ is {\it unipotent} with respect to $\varphi$
if $\partial{(x \circ \varphi)}/\partial{x} \equiv 1$ in $\gamma$.
If $\gamma$ is unipotent then the germ of $\varphi$ at $P$ is unipotent
for all $P \in \gamma$. We will prove that $\log \varphi$ is u.t.f. along $\gamma$.

Let $P \in Fix \varphi$. Let $\varphi_{P}$ be the 1-dimensional germ
obtained in the neighborhood of $P$ by restricting $\varphi$ to $\cap_{j=1}^{n}(x_{j}=x_{j}(P))$.
Then $\varphi_{P}$ can be embedded in a formal flow
except if $|(\partial{(x \circ \varphi)}/\partial{x})(P)|$ is a root of the unit
different than $1$ and $\varphi_{P}$ is not periodic. As a consequence there is no
hope in general for $\log \varphi$ to be u.t.f. along the non-unipotent components of $Fix \varphi$.
Nevertheless, if ${\varphi}_{P}$
can be embedded in a formal flow for all $P \in \gamma$ contained in  an irreducible component $\gamma$
of $Fix \varphi$ then $\log \varphi$ is u.t.f. along $\gamma$.
We will make clear the previous assertions.
\subsection{One-dimensional results}
The next results are well-known and they are included here for the sake of clarity.
\begin{pro}
\label{pro:explin}
  Let $\tau \in \diff{}{}$ such that $j^{1} \tau \neq Id$. Then $\tau$ is the exponential of a
formal vector field if and only if it is formally linearizable.
\end{pro}
\begin{proof}
Let $\hat{X}=(a_{1} z + a_{2} {z}^{2} + \hdots) \partial / \partial{z}$ be a formal vector field
such that $\tau = {\rm exp}(\hat{X})$. Since $j^{1} \tau = {e}^{a_{1}} z$ and $j^{1} \tau \neq Id$
then $a_{1} \neq 0$.
The linear part is a complete system of invariants for the elements of $\hat{\mathcal X} \cn{}$
with no vanishing linear part. As a consequence there exists $\hat{\sigma} \in \diffh{}{}$ such that
$\hat{\sigma}_{*} (a_{1} z \partial/ \partial z)= \hat{X}$. By taking exponentials we obtain
$\hat{\sigma} \circ j^{1} \tau = \tau \circ \hat{\sigma}$.

Suppose $\hat{\sigma} \circ j^{1} \tau = \tau \circ \hat{\sigma}$. Since
$j^{1} \tau = {\rm exp}(c z \partial / \partial{z})$ for all $c \in {\mathbb C}$
such that ${e}^{c}= (\partial{\tau} / \partial{z})(0)$
then $\tau = {\rm exp}(\hat{\sigma}_{*}(c z \partial / \partial{z}))$.
\end{proof}
\begin{pro}
\label{pro:unqlog}
Let $\tau \in \diff{}{}$ with $j^{1} \tau$ not periodic. Assume that we have
$\tau = {\rm exp}(\hat{X}) = {\rm exp}(\hat{Y})$ for $\hat{X}$, $\hat{Y} \in \hat{\mathcal X} \cn{}$ such
that $j^{1} \hat{X} = j^{1} \hat{Y}$. Then $\hat{X} = \hat{Y}$.
\end{pro}
\begin{proof}
We can suppose $\tau = j^{1} \tau$ by proposition \ref{pro:explin}.
It is enough to prove that $\hat{X} = j^{1} \hat{X}$.
Since $j^{1} \hat{X} \neq 0$ there exists $\hat{\sigma} \in \diffh{}{}$ such that
$\hat{\sigma}_{*} (j^{1} \hat{X}) = \hat{X}$.
By taking exponentials we obtain $\hat{\sigma} \circ \tau = \tau \circ \hat{\sigma}$.
Since $\tau$ is linear and not periodic the center of $\tau$ in $\diffh{}{}$ is the
linear group. Therefore $\hat{\sigma}$ is linear, that implies
$j^{1} \hat{X} = \hat{\sigma}_{*} (j^{1} \hat{X}) = \hat{X}$.
\end{proof}
\begin{pro}
\label{pro:concut}
Let $\tau \in \diff{}{}$ such that $j^{1} \tau$ is periodic and $\tau$ is linearizable.
Suppose we have ${h} \circ \tau \circ h^{(-1)}(z) - j^{1} \tau(z) = O({z}^{k+1})$ for some
$h \in \diff{}{}$ and some $k \geq 1$. Then there exists $\sigma \in \diff{}{}$ such that
$j^{k} \sigma = j^{k} h$ and $\sigma \circ \tau \circ {\sigma}^{(-1)} = j^{1} \tau$.
\end{pro}
\begin{proof}
It is enough to prove that if $\tau - j^{1} \tau= O({z}^{k+1})$ there exists
$\sigma \in \diff{}{}$ such that $\sigma \circ \tau = j^{1} \tau \circ \sigma$
and $\sigma(z) - z =O({z}^{k+1})$. We denote $j^{1} \tau = a z$ and
the period of $j^{1} \tau$ by $q$. Since
$\tau^{(q)}$ is formally conjugated to ${(j^{1} \tau)}^{(q)} = Id$ then $\tau^{(q)}=Id$.
We are done by defining $\sigma(z) = (\sum_{j=0}^{q-1} \tau^{(j)}(z)/a^{j})/q$.
\end{proof}
\subsection{Division neighborhoods. Convergence by restriction}
We work in domains in which the components of $Fix \varphi$ behave like
their germs at $0$.
Consider $g \in {\mathbb C} \{ x, x_{1}, \hdots , x_{n} \}$. We say that a domain $U$ is a
{\it division neighborhood} for $g$ if
\begin{itemize}
\item There is a decomposition $g = g_{1}^{l_{1}} \hdots g_{r}^{l_{r}}$ of $g$ in irreducible
factors in the ring ${\mathbb C} \{ x, x_{1}, \hdots , x_{n} \}$ such that
$g_{j} \in \vartheta(U)$ for $1 \leq j \leq r$. \\
\item The regular part of $g_{j}=0$ is connected in $U$ for all $1 \leq j \leq r$.
\end{itemize}
The definition of division neighborhood is intended to extend the division of germs to
bigger domains. More precisely, if $U$ is a division neighborhood for $g$ then
\[ a = g b \ {\rm where} \ a \in {\vartheta}(U) \ {\rm and} \
b \in {\mathbb C}\{x,x_{1},\hdots,x_{n}\} \implies b \in {\vartheta}(U) . \]
The following results can be immediately deduced from the definition of division neighborhood.
\begin{lem}
\label{lem:sub}
Let $U$ be a division neighborhood for $g$ and consider a domain $V \subset U$
such that $U \cap (g=0) = V \cap (g=0)$. Then $V$ is a division neighborhood for $g$.
\end{lem}
\begin{lem}
\label{lem:remana}
Let $U$ be a division neighborhood for $g$. Consider an analytic set $S \not \ni 0$.
Then $U \setminus S$ is a division neighborhood for $g$.
\end{lem}
The last lemma is trivial since we can not break the connectedness of the
regular parts by removing sets of real codimension at least $2$.

It is not difficult to find
a division neighborhood for $f=x \circ \varphi -x$.
Let $(y_{1},\hdots , y_{n+1})$ be a set of coordinates such that any irreducible component
$\gamma$ of $Fix \varphi$ can be expressed as the vanishing set of a monic Weierstrass polynomial
in the variable $y_{1}$. Let
$(c_{1},\hdots,c_{n+1}) \in {\mathbb R}^{+} \times \hdots \times {\mathbb R}^{+}$.
Every polydisk $\cap_{j=1}^{n+1} (|y_{j}| < c_{j})$ small enough
such that
\[ (|y_{1}|=c_{1}) \cap \cap_{j=2}^{n+1} (|y_{j}| < c_{j}) \cap (f=0) = \emptyset \]
is a division neighborhood for $f$.

  We define a new concept, it is simpler to handle than the formal transversality; it is a sort
of formal transversality at the $0$-level. We say that
$h \in {\mathbb C} [[ x, x_{1}, \hdots , x_{n} ]]$
{\it converges by restriction to}
$\gamma$ if there exists $h' \in {\mathbb C} \{ x, x_{1}, \hdots , x_{n} \}$
such that $h - h'$ belongs to $I(\gamma)$ where $I(\gamma)$ is the ideal of $\gamma$.
If $V \subset \gamma$ is  a neighborhood of $0$ in $\gamma$
we say that $h$ {\it converges by restriction to} $\gamma$ {\it in} $V$ if
$h'$ can be chosen holomorphic in a neighborhood of $V$.
\begin{lem}
\label{lem:cbr}
  Let $h$ be an irreducible element of ${\mathbb C} \{ x, x_{1}, \hdots , x_{n} \}$.
Consider series $\hat{c} \in {\mathbb C} [[ x, x_{1}, \hdots , x_{n} ]]$ and
$d \in {\mathbb C} \{ x, x_{1}, \hdots , x_{n} \} \setminus (h)$. Suppose $\hat{c} d$ converges
by restriction to $h=0$. Then $\hat{c}$ converges by restriction to $h=0$.
\end{lem}
There is also a uniform version of the previous lemma but we have to be careful with the
setting. We keep the notations in the previous lemma. We consider coordinates
$(y_{1},\hdots,y_{n+1})$ such that
\[ h = v(y_{1}, \hdots , y_{n+1}) (y_{1}^{l} + y_{1}^{l-1} a_{l-1}(y_{2},\hdots,y_{n+1}) +
\hdots +  a_{0}(y_{2},\hdots,y_{n+1})) \]
for some unit $v$ and some functions $a_{j}$ ($0 \leq j \leq l-1$). For a generic
point $(y_{2}^{0}, \hdots , y_{n}^{0})$ the number of points in
$h(y_{1}, y_{2}^{0}, \hdots , y_{n}^{0}) = 0$ is $l$. We enumerate them
${\alpha}_{1}(y_{2},\hdots,y_{n+1})$, $\hdots$, ${\alpha}_{l}(y_{2},\hdots,y_{n+1})$.
We define
\[ \Delta(h,d)  = \left({ \prod_{j=1}^{l} (d \circ {\alpha}_{j})
\prod_{1 \leq j < k \leq l} {({\alpha}_{j}  - {\alpha}_{k} ) }^{2} }\right) (y_{2}, \hdots , y_{n+1}).  \]
The function $\Delta(h,d)$ is well-defined because it is symmetric in the functions
${\alpha}_{1}$, $\hdots$, ${\alpha}_{l}$; it is holomorphic in $\Delta(h,d) \neq 0$
and continuous in a neighborhood of the origin. By Riemman's theorem
$\Delta(h,d)$ is holomorphic in a neighborhood of the origin.
There exists a division neighborhood $D= D_{1} \times D_{2, \hdots, n+1} \subset {\mathbb C} \times {\mathbb C}^{n}$
for $h$ such that $D_{2, \hdots, n+1}$ is a division neighborhood
for $\Delta$. Given any neighborhood of the origin
$W$ the set $D$ can be chosen to be contained in $W$. Suppose that $d$ converges in a neighborhood of
$D \cap (h=0)$. In this context we prove
\begin{lem}
\label{lem:cbru}
  If $\hat{c} d$ converges by restriction to $h=0$ in $D \cap (h=0)$ and $d \not \in (h)$
then $\hat{c}$ converges by restriction to $h=0$ in $D \cap (h=0)$.
\end{lem}
\begin{proof}[proof of lemmas \ref{lem:cbr} and \ref{lem:cbru}]
It is clear that the lemma \ref{lem:cbru} implies the lemma \ref{lem:cbr}.
  Consider a holomorphic function $c_{0}$ defined in the neighborhood of $D \cap (h=0)$
such that $\hat{c} d - c_{0} \in (h)$. Next step is
using the Weierstrass division, we want to divide $c_{0} / d$ by $h$. The remainder of
that division is
\[ R = \sum_{j=1}^{l} \frac{c_{0}(\alpha_{j}(y_{2},\hdots,y_{n+1}))}{d(\alpha_{j}(y_{2},\hdots,y_{n+1}))}
\frac{\prod_{k \neq j} (y_{1} - {\alpha}_{k}(y_{2},\hdots,y_{n+1}))}
{\prod_{k \neq j} ({\alpha}_{j}(y_{2},\hdots,y_{n+1}) - {\alpha}_{k}(y_{2},\hdots,y_{n+1}))} . \]
The function $\Delta R$ is holomorphic in a neighborhood of the origin
and since $\Delta$ does not depend on $y_{1}$ then
$\Delta R$ is the remainder of the Weierstrass division $[(\Delta c_{0}) / d]/h$.
We define $\hat{R}$ the remainder of the Weierstrass division $\hat{c}/h$, it is a formal
power series. We have
\[ \Delta c_{0} - \Delta R d \in (h) \ \ {\rm and} \ \ \Delta \hat{c} d - \Delta \hat{R} d \in (h) . \]
Since $d \not \in (h)$ we obtain $\Delta R - \Delta \hat{R} \in (h)$. Both $\Delta R$ and
$\Delta \hat{R}$ are polynomials in the variable $y_{1}$ whose degree is lesser or equal than $l-1$.
Therefore, we have $\Delta R \equiv \Delta \hat{R}$ by uniqueness of the Weierstrass division.
The function $\Delta$ divides the $l$ coefficients of the polynomial $\Delta R$. As a consequence
$R$ is convergent in the neighborhood of the origin. Since
$D_{2,\hdots,n+1}$ is a division neighborhood for $\Delta$ then $R$ is defined
in ${\mathbb C} \times D_{2,\hdots,n+1}$. The series $\hat{c}$ converges by restriction to $h=0$
in $D \cap (h=0)$ since $\hat{c} - R = \hat{c} - \hat{R}   \in (h)$.
\end{proof}
We relate formal transversality along an analytic hypersurface with formal transversality along its
irreducible components.
\begin{lem}
\label{lem:allforone}
  Let $f$ be an element of ${\mathbb C} \{ x, x_{1}, \hdots , x_{n} \}$. Then
$\hat{u}$ in ${\mathbb C} [[ x, x_{1}, \hdots , x_{n} ]]$ is t.f. (resp. u.t.f.) along $f=0$
if and only if $\hat{u}$ is t.f. (resp. u.t.f.) along every irreducible component of $f=0$.
\end{lem}
\begin{proof}
The sufficient condition is obvious.

Suppose we are in the u.t.f. case, the proof for the t.f. case is simpler.
  Let $f=f_{1}^{n_{1}} \hdots f_{l}^{n_{l}}$ be the decomposition of $f$ in irreducible factors
in ${\mathbb C} \{ x, x_{1}, \hdots , x_{n} \}$.
We choose coordinates $(y_{1}, \hdots, y_{n+1})$ such that
$f_{j}=0$ can be expressed up to a unit as a monic Weierstrass polynomial in the variable $y_{1}$
for all $1 \leq j \leq l$. We can choose a polydisk $D= D_{1} \times D_{2,\hdots, n+1}$
such that $(\partial{D_{1}} \times D_{2,\hdots,n+1}) \cap (f=0) = \emptyset$.
Moreover, we can suppose that
$D$ is a division neighborhood for every $f_{j}$ and $D_{2,\hdots,n+1}$ is
a division neighborhood for every $\Delta(f_{j},f_{k})$ with $j \neq k$. Finally we
suppose that $\hat{u} \in \lim_{\leftarrow} \vartheta(D)/(f_{j})^{k}$
for all $1 \leq j \leq l$.

Let $F=\prod_{k=1}^{l} f_{k}^{a_{k}}$. Suppose we have $u \in \vartheta(D)$
such that $\hat{u} -u \in (F)$. Fix $j \in \{1, \hdots, l\}$; it is enough to prove the
existence of a function $v \in  \vartheta(D)$ such that
$\hat{u} -v \in (f_{j} F)$.
We claim that $(\hat{u} -u) / f_{j}^{a_{j}}$ converges by restriction to $f_{j}=0$.
There exists $u_{a_{j}} \in \vartheta(D)$ such
that $\hat{u}-u_{a_{j}}$ belongs to $({f_{j}}^{a_{j}+1})$ by hypothesis.
Denote $w= (u_{a_{j}} - u) / f_{j}^{a_{j}}$; we have
$(\hat{u} - u) / f_{j}^{a_{j}} - w \in (f_{j})$.
By lemma \ref{lem:cbru} the series $(\hat{u}-u)/F$ converges by restriction to
$f_{j}=0$ in $D \cap (f_{j}=0)$. As a consequence there exists
$b(y_{1},\hdots,y_{n+1}) \in \vartheta({\mathbb C} \times D_{2,\hdots,n+1}) \subset \vartheta(D)$
such that $(\hat{u}-u)/F - b \in (f_{j})$. This implies $\hat{u} -(u + bF) \in (f_{j}F)$.
\end{proof}
Consider a set $W \subset {\mathbb C}^{n+1}$. We can define the ring $G_{W}$ of germs of functions
in a neighborhood of $W$. The next lemma provides a handy characterization of u.t.f. functions.
%
%
%
%
\begin{lem}
  Let $V \subset {\mathbb C}^{n+1}$ be a germ of analytic hypersurface at $0$. Then a series
$\hat{g} \in {\mathbb C}[[x,x_{1},\hdots,x_{n}]]$ is u.t.f. along $V$
if and only if $\hat{u}$ belongs to $\lim_{\leftarrow} G_{V \cap W}/I(V)^{j}$
for some neighborhood $W$ of the origin.
%
%
%
%
%
\end{lem}
\begin{proof}
The sufficient condition is obvious.

Consider coordinates $(y_{1}, \hdots , y_{n+1})$ such that $I(V)=(h)$ for some
monic Weierstrass polynomial $h \in {\mathbb C}[y_{1}][[y_{2},\hdots,y_{n+1}]]$ in the variable $y_{1}$.
Choose a polydisk $D=D_{1} \times D_{2,\hdots,n+1} \subset {\mathbb C} \times {\mathbb C}^{n}$
in the variables $(y_{1},\hdots,y_{n+1})$ such that $(\partial{D_{1}} \times D_{2,\hdots,n+1}) \cap V =\emptyset$.
Moreover, we choose $D$ such that it is a division neighborhood for $h$ and $D \cap V \subset W$. We denote
$\pi(y_{1},\hdots , y_{n+1}) = (y_{2}, \hdots , y_{n+1})$.
There exists a function $b_{j} \in G_{V \cap W}$ such that $\hat{g} - b_{j} \in (h^{j})$ for all
$j \in {\mathbb N}$. Since $h^{j}$ is a monic polynomial in $y_{1}$
we can consider the remainder $g_{j}$ of the Weierstrass division $b_{j}/h^{j}$.
Since $b_{j} \in G_{D \cap V}$
then $g_{j} \in \vartheta({\pi}^{-1} (\pi (D \cap V))) \subset \vartheta(D)$
for all $j \in {\mathbb N}$. Clearly $\hat{g}$ is u.t.f. along $V$.
%
%
%
%
\end{proof}
\begin{rem}
  The results in this section can be enounced in the uniform semi-meromorphic case with minor
adjustments. For the sake of simplicity we omit such a formulation.
\end{rem}
\subsection{Main results}
Fix $\varphi = {\rm exp}(\hat{u} (x \circ \varphi -x) \partial/\partial{x}) \in \diff{up}{n+1}$
and an irreducible component $\gamma$ of $Fix \varphi$.
Consider an irreducible equation $g=0$ of $\gamma$.
Denote $f=x \circ \varphi - x$.
Let $e$ be the multiplicity of
$\gamma$ in $f=0$, i.e. the greatest $e \in {\mathbb N}$ such that $g^{e}$ divides $f$.
Denote $f/{g}^{e}$ by $h$.
We choose a set of coordinates $(y_{1},\hdots , y_{n+1})$ and a polydisk
$D = D_{1} \times D_{2, \hdots , n+1} \subset {\mathbb C} \times {\mathbb C}^{n}$ in
the variables $(y_{1},\hdots , y_{n+1})$ such that
$(\partial{D_{1}} \times D_{2,\hdots,n+1}) \cap \gamma=\emptyset$ and
\begin{itemize}
\item $\varphi$ is holomorphic in $D$. \\
\item Up to a multiplicative unit $g$ is of the form
\[ g = y_{1}^{l} + y_{1}^{l-1} a_{l-1}(y_{2},\hdots,y_{n+1}) + \hdots + a_{0}(y_{2},\hdots,y_{n+1}) \]
where $a_{j}(0)=0$ for all $0 \leq j \leq l -1$. \\
\item $D$ is a division neighborhood for $f=0$. \\
\item $D_{2,\hdots,n+1}$ is a division neighborhood for $\Delta(g,h)$.
\end{itemize}If $\gamma$ is a non-unipotent component there are two more conditions:
\begin{itemize}
\item $D_{2,\hdots,n+1}$ is a division neighborhood for $\Delta(g, \partial{f}/\partial{x})$. \\
\item $\ln( \partial(x \circ \varphi) / \partial{x})$ is a holomorphic function defined in $D$
such that the set $[\partial (x \circ \varphi) / \partial{x} =1] \cap (D \cap \gamma)$ is connected.
\end{itemize}
These notations are fixed throughout this section.
\begin{pro}
\label{pro:tfuhea}
  Let $\varphi$ be a up-diffeomorphism and let $\gamma$ be a unipotent irreducible component of
$Fix \varphi$. Then $\log \varphi$ is u.t.f. along $\gamma$.
\end{pro}
\begin{proof}
The proprieties $e>1$ if
$\gamma$ is non-fibered or $g | \partial{g}/\partial{x}$ if $\gamma$ is fibered imply the
existence of a formal series $\hat{w}$ such that
\[ \varphi = (x + \hat{u} {g}^{e} h + \hat{w} g^{e+1} h,x_{1}, \hdots , x_{n+1}) . \]
It is a consequence of the Taylor series expansion for the exponential mapping.
Hence $u_{1} = 1=(x \circ \varphi - x)/f \in \vartheta(D)$ satisfies $\hat{u} - u_{1} \in (g)$.

  We proceed by induction. Suppose $\hat{u} = u_{\alpha} + {g}^{\alpha} \hat{v}$ where
$u_{\alpha} \in G_{D \cap \gamma}$; the result is already
proved for $\alpha=1$. We want to find a similar expression for $\alpha+1$. It is enough
to prove the existence of $v \in G_{D \cap \gamma}$ such that
$\hat{v} - v \in (g)$. We have
\[ x \circ \varphi = x \circ {\rm exp} (u_{\alpha} f \partial{}/\partial{x}) +
{g}^{\alpha + e} h \hat{v} + {g}^{\alpha + e +1} h \hat{A}  \]
for a certain formal series $\hat{A}$. The vector field $u_{\alpha} f \partial{}/\partial{x}$
is defined in a neighborhood of $D \cap \gamma$ and it vanishes in $\gamma$.
Therefore, its exponential is defined in a neighborhood of $D \cap \gamma$.
We define
\[ d = (x \circ \varphi - x \circ {\rm exp} (u_{\alpha} f \partial{}/\partial{x}))/ {g}^{\alpha +e}. \]
We obtain that $d \in G_{D \cap \gamma}$ by lemma \ref{lem:sub} since $D$ is a division
neighborhood for $f$ and then for $g$. We have $h \hat{v} - d \in (g)$, thus $h \hat{v}$
converges by restriction to $\gamma$
in $D \cap \gamma$. By applying lemma \ref{lem:cbru} we get that $\hat{v}$
converges by restriction to $\gamma$ in $D \cap \gamma$. Then there exists
$v \in G_{D \cap \gamma}$ such that $\hat{v} - v \in (g)$.
We define $u_{\alpha +1} = u_{\alpha} + g^{\alpha} v$, we have $\hat{u} - u_{\alpha +1} \in ({g}^{\alpha +1})$.
\end{proof}
Fix $\varphi = {\rm exp}(\hat{u} (x \circ \varphi-x) \partial / \partial{x}) \in \diff{up}{n+1}$ and an irreducible
component $\gamma$ of $Fix \varphi$. We can express $\hat{u}$ in the form
$\sum_{j=0}^{\infty} u_{j}^{\gamma}  {g}^{j}$
where $u_{j}^{\gamma} \in {\mathbb C}[y_{1}][[y_{2},\hdots,y_{n+1}]]$ is a polynomial
in $y_{1}$ such that $\deg_{y_{1}} u_{j}^{\gamma} \leq l-1$
for all $j \geq 0$. These properties characterize the sequence $\{ u_{j}^{\gamma} \}$
since $u_{0}^{\gamma}$ is the remainder of the Weierstrass division $\hat{u} / g$,
the series $u_{1}^{\gamma}$ is the Weierstrass remainder of $[(\hat{u}-u_{0}^{\gamma})/g]/g$ and so on.

We define $E_{0}^{\gamma}=E_{1}^{\gamma}= \emptyset$ and
\[ E_{k}^{\gamma} = \left\{{ Q \in \gamma \ : \ \exists \ 1 < j \leq k \ {\rm s.t.}  \
{\left({ \frac{\partial{(x \circ \varphi)}}{\partial{x}} }\right)}^{j} (Q)  = 1   }\right\}
\setminus
\left\{{ \frac{\partial{(x \circ \varphi)}}{\partial{x}} = 1}\right\} \]
for all $k \geq 2$. We define $\pi(y_{1}, \hdots, y_{n+1})=(y_{2}, \hdots, y_{n+1})$.
\begin{lem}
\label{lem:genexp}
(Generic expression)
Let $\varphi = {\rm exp}(\hat{u} f \partial/\partial{x})$ be a up-diffeomorphism
and let $\gamma$ be a non-unipotent irreducible component of
$Fix \varphi$.
Then $u_{j}^{\gamma}$ is holomorphic in
${\pi}^{-1}(D_{2,\hdots,n+1} \setminus \pi(E_{j}^{\gamma}))$ for all $j \geq 0$.
In particular $u_{j}^{\gamma} \in {\mathbb C} \{ x, x_{1}, \hdots , x_{n} \}$ for all $j \geq 0$.
Moreover, we have
\[ \varphi = {\rm exp} \left({
(\sum_{j=0}^{\infty} u_{j}^{\gamma} {g}^{j}) g h \frac{\partial{}}{\partial{x}}
}\right)  \]
in the neighborhood of every point in
$\gamma \cap [D \setminus \pi^{-1}(\overline{\cup_{j \geq 0} \pi(E_{j}^{\gamma})})]$.
\end{lem}
\begin{proof}
We have
\begin{equation}
\label{equ:genexp}
 x \circ \varphi - \left({
x + \frac{1}{k!} \sum_{k=1}^{\infty} {(u_{0}^{\gamma})}^{k} h^{k} {(\partial{g}/\partial{x})}^{k-1} g }\right)
\in (g^{2})
\end{equation}
and then $\partial (x \circ \varphi) / \partial{x} -
{e}^{u_{0}^{\gamma} h \partial{g}/\partial{x}} \in (g)$. We deduce that
$u_{0}^{\gamma} h \partial{g}/\partial{x}$ converges by restriction to $\gamma$ in $D \cap \gamma$.
By lemma  \ref{lem:cbru} there exists $u_{0}' \in G_{D \cap \gamma}$
such that $u_{0}^{\gamma} - u_{0}' \in (g)$. Since $u_{0}^{\gamma}$
coincides with the remainder of the Weierstrass division $u_{0}'/g$ then
$u_{0}^{\gamma}$ is holomorphic in ${\pi}^{-1}(D_{2,\hdots,n+1})$.
We denote $L= {e}^{u_{0}^{\gamma} h \partial{g}/\partial{x}}$.

  We are going to prove the result by induction on $j$. The result is true for $u_{0}^{\gamma}$.
Suppose $u_{j}^{\gamma}$ is holomorphic in ${\pi}^{-1}(D_{2,\hdots,n+1} \setminus \pi(E_{j}^{\gamma}))$
for all $j < \alpha$. We will prove that
$u_{\alpha}^{\gamma}$ is holomorphic in ${\pi}^{-1}(D_{2,\hdots,n+1} \setminus \pi(E_{\alpha}^{\gamma}))$.
Let $u=  \sum_{j=0}^{\alpha -1} u_{j}^{\gamma} g^{j}$; we denote
$u_{0}=u_{0}^{\gamma}$ and $\hat{v} =(\hat{u}-u)/{g}^{\alpha}$.

  Let $\hat{B}_{0}=H_{0}=x$; we define
$\hat{B}_{j+1} = \hat{u} g h \partial{\hat{B}_{j}} / \partial{x}$ and
$H_{j+1}= u g h \partial{H_{j}} / \partial{x}$ for $j \geq 0$. The next step in the
proof is proving by induction that $\hat{B}_{j}$ ($j \geq 0$) can be expressed in the form
\[ \hat{B}_{j} = H_{j} + \hat{v}  {g}^{\alpha +1} h^{j} u_{0}^{j-1}
{(\partial{g} / \partial{x})}^{j-1} C_{j} + {g}^{\alpha +2} h \hat{D}_{j} \]
where $\hat{D}_{j}$ belongs to ${\mathbb C} [[ x, x_{1}, \hdots , x_{n} ]]$ and
$C_{j} \in {\mathbb C}$. The relations defining the sequence
are $C_{0}=0$ and $C_{j}=(\alpha + 1) C_{j-1} + 1$ for all $j > 0$.
The induction result is immediate for $j=0$ and $j=1$. We can develop
\[ \hat{B}_{j+1} = (\partial{\hat{B}_{j}} / \partial{x}) ( u + {g}^{\alpha} \hat{v}) g h \]
to obtain
\[ \hat{B}_{j+1} \equiv \left[{
\frac{\partial{H_{j}}}{\partial{x}} + \frac{\partial{}}{\partial{x}}
\left({
\hat{v} {g}^{\alpha +1} h^{j} u_{0}^{j-1} {\left({ \frac{\partial{g}}{\partial{x}} }\right)}^{j-1} C_{j}
}\right)
}\right] (u g h + {g}^{\alpha + 1} h \hat{v}) \]
modulo $({g}^{\alpha +2} h)$. By using $u - u_{0} \in (g)$ and the equation \ref{equ:genexp} we get
\[ \partial{H_{j}} / \partial{x} - u_{0}^{j} h^{j} {(\partial{g}/\partial{x})}^{j} \in (g) . \]
This leads us to the desired expression
\[ \hat{B}_{j+1} - [H_{j+1} + u_{0}^{j} h^{j+1} {(\partial{g}/\partial{x})}^{j} {g}^{\alpha +1}
\hat{v} (1 + C_{j}(\alpha + 1))] \in ({g}^{\alpha +2} h) . \]
for $\hat{B}_{j+1}$. We note that $C_{j}= ({(\alpha + 1 )}^{j} -1)/ \alpha$.
We have
\[ x \circ \varphi \equiv x \circ {\rm exp}(u g h \partial{}/\partial{x}) +
\left({ \sum_{j=1}^{\infty} \frac{C_{j}  u_{0}^{j-1} h^{j-1} {(\partial{g}/\partial{x})}^{j-1}}{j!}
}\right)  h g^{\alpha +1} \hat{v}  \]
modulo $(g^{\alpha +2} h)$. We simplify to obtain
\[
x \circ \varphi - \left[{
x \circ {\rm exp}(u g h \partial{}/\partial{x}) +
\frac{(L^{\alpha +1} - L) h g^{\alpha +1} \hat{v}}{\alpha u_{0} h \partial{g} / \partial{x}}
}\right] \in ({g}^{\alpha +2} h) .
\]
Since $u$ is defined in a neighborhood of
$(D \cap \gamma) \setminus {\pi}^{-1}(\pi(E_{\alpha -1}^{\gamma}))$ and $u g h \partial{}/\partial{x}$
vanishes on $\gamma$ then $x \circ {\rm exp}(u g h \partial{}/\partial{x})$ is defined
in a neighborhood of $(D \cap \gamma) \setminus \pi^{-1} (\pi(E_{\alpha -1}^{\gamma}))$.
We have the inclusion
\[ [\partial (x \circ \varphi) / \partial{x} =1] \cap D \cap \gamma =
[L=1] \cap D \cap \gamma \subset [u_{0} h \partial{g} / \partial{x} = 0] \]
since otherwise $[\partial (x \circ \varphi) / \partial{x} =1] \cap D \cap \gamma$ is not connected.
As a consequence we obtain that $(L-1)/(u_{0} h \partial{g} / \partial{x})$ is never vanishing
in $D \cap \gamma$. The function
\[ K \stackrel{def}{=}
\frac{x \circ \varphi - x \circ {\rm exp}(u g h \partial{}/\partial{x})}{1 + L + \hdots + L^{\alpha - 1}}
\frac{\alpha}{L} \frac{u_{0} h \partial{g} / \partial{x}}{L-1}  \]
is holomorphic in a neighborhood of
$(D \cap \gamma) \setminus (\pi^{-1} (\pi(E_{\alpha -1}^{\gamma})) \cup E_{\alpha}^{\gamma})$.
The polydisk $D$ is a division neighborhood for $g$. By lemmas \ref{lem:remana} and \ref{lem:sub}
the series $K' = K / {g}^{\alpha +1}$ is a holomorphic function in a neighborhood of
$(D \cap \gamma) \setminus (\pi^{-1} (\pi(E_{\alpha -1}^{\gamma})) \cup E_{\alpha}^{\gamma})$.

Since $(K' / h) h$ converges by restriction to $\gamma$ in
$(D \cap \gamma) \setminus \pi^{-1} (\pi(E_{\alpha}^{\gamma}))$ then so does $K'/h$ by lemmas
\ref{lem:cbru} and \ref{lem:remana}. The series $u_{\alpha}^{\gamma}$ is the Weierstrass
remainder of the division $[K' / h] / g$. Therefore $u_{\alpha}^{\gamma}$ is holomorphic
in ${\pi}^{-1}(D_{2,\hdots,n+1} \setminus \pi(E_{\alpha}^{\gamma}))$.
%
%

Since $x \circ \varphi - x \circ {\rm exp}( (\sum_{j=0}^{\alpha-1} u_{j}^{\gamma}
{g}^{\alpha-1})g h \partial/\partial{x}) \in ({g}^{\alpha +1})$ we can apply
lemmas \ref{lem:remana} and \ref{lem:sub}
to prove the existence of a holomorphic function $M_{\alpha}$ defined in a neighborhood of
$(D \cap \gamma) \setminus \pi^{-1} (\pi(E_{\alpha-1}^{\gamma}))$ and such that
\[ x \circ \varphi -
x \circ {\rm exp} \left({ (\sum_{j=0}^{\alpha-1} u_{j}^{\gamma} {g}^{\alpha-1})g h \partial/\partial{x} }\right)
= {g}^{\alpha+1} M_{\alpha} . \]
The previous equality running on $\alpha \in {\mathbb N}$ imply that
\[ \varphi = {\rm exp} \left({
(\sum_{j=0}^{\infty} u_{j}^{\gamma} {g}^{j}) g h \frac{\partial{}}{\partial{x}}
}\right)  \]
in the neighborhood of every point in
$\gamma \cap [D \setminus \pi^{-1}(\overline{\cup_{j \geq 0} \pi(E_{j}^{\gamma})})]$.
\end{proof}
\begin{rem}
  There is not always uniform formal transversality. The logarithm of
$\varphi = (x + y - {x}^{2} , y)$ is not u.t.f. along $\gamma \equiv (y={x}^{2})$.
We claim that ${\varphi}_{P}$
can not be embedded in a formal flow if $P \in \cup_{j \geq 0} E_{j}^{\gamma}$.
Otherwise $\varphi_{P}$ is formally conjugated to $j^{1} \varphi_{P}$ (prop. \ref{pro:explin})
and then periodic; that is not possible since ${\varphi}_{P}$ is non-linear and polynomial.
We obtain that $\log \varphi$ is not u.t.f. along $\gamma$ since $0 \in \overline{\cup_{j \geq 0} E_{j}^{\gamma}}$.
\end{rem}
\begin{pro}
\label{pro:pretra}
Let $\varphi = {\rm exp}(\hat{u} f \partial/\partial{x}) \in \diff{up}{n+1}$
and let $\gamma$ be a non-unipotent irreducible component of
$Fix \varphi$. The coefficients of
$u_{j}^{\gamma}$ in ${\mathbb C} \{y_{2}, \hdots y_{n+1}\}$ are meromorphic in
$D_{2,\hdots,n+1}$ for all $j \geq 0$.
\end{pro}
\begin{proof}
We denote $L={e}^{u_{0}^{\gamma} h \partial{g}/\partial{x}}$; we know that
$L - \partial{(x \circ \varphi)}/\partial{x} \in (g)$.
Let $P$ be a point of $(D \cap \gamma) \setminus (L=1)$. Since
$(L=1) \cap \gamma = (\partial{f}/\partial{x}=0) \cap \gamma$ then $f=0$ is transversal
to $\partial{} / \partial{x}$ in $P$. In particular $\gamma$ is smooth at $P$ and it
is the only irreducible component of $Fix \varphi$ passing through $P$.
We can find
new coordinates $(z,x_{1},\hdots,x_{n})$ in the neighborhood of $P$ such that
$\partial{}/\partial{x}$ and $g=0$ become $\partial{}/\partial{z}$ and $z=0$ respectively.
In these coordinates $\varphi$ is of the form
\[ \varphi(z,x_{1},\hdots,x_{n})= (a_{1}^{1}(x_{1}, \hdots , x_{n})z + a_{2}^{1}(x_{1}, \hdots , x_{n}) {z}^{2}
+ \hdots ,x_{1},\hdots,x_{n}) . \] Consider a small enough open neighborhood
$W(P)$ of $P$ in $D \cap \gamma$, it is parameterized by $(x_{1}, \hdots , x_{n})$.
We ask $a_{1}^{1} - 1$ for never vanishing
in $W(P)$; this is possible since $L(P) \neq 1$.
We claim that there exists $\sigma$  such that
$\sigma^{-1} \circ \varphi \circ \sigma = (a_{1}^{1} z, z_{1}, \hdots , z_{n})$
of the form
\[ \sigma = \left({ z + \sum_{j=2}^{\infty} \sigma_{j}(x_{1}, \hdots x_{n}) {z}^{j}, x_{1}, \hdots , x_{n} }\right), \]
where $\sigma_{j}$ is holomorphic in $W(P) \setminus E_{j-1}^{\gamma}$
and meromorphic in $W(P)$ for all $j \geq 2$.
We are going to construct a sequence of diffeomorphisms $(\varphi_{k})$
such that $\varphi_{1}=\varphi$ and $\varphi_{k}$ is of the form
\[ \varphi_{k} = \left({ a_{1}^{1}(x_{1}, \hdots , x_{n})z +
\sum_{j=k+1}^{\infty} a_{j}^{k}(x_{1}, \hdots , x_{n}) {z}^{j}, x_{1}, \hdots , x_{n} }\right) \]
where $a_{j}^{k}$ is holomorphic in $W(P) \setminus  E_{k-1}^{\gamma}$
and meromorphic in $W(P)$ for all $j \geq 2$. There
exists $\tau_{k}=(z+b_{k+1}(x_{1},\hdots,x_{n}) {z}^{k+1},x_{1}, \hdots , x_{n})$ such that
\[ z \circ \tau_{k}^{(-1)} \circ \varphi_{k} \circ \tau_{k} - a_{1}^{1}(x_{1}, \hdots , x_{n})z \in ({z}^{k+2}) . \]
Moreover $\tau_{k}$ is unique and $({(a_{1}^{1})}^{k+1} - a_{1}^{1}) b_{k+1}=a_{k+1}^{k}$.
By hypothesis $b_{k+1}$ is holomorphic in $W(P) \setminus E_{k}^{\gamma}$
and meromorphic in $W(P)$. We define the diffeomorphism
$\varphi_{k+1}= \tau_{k}^{(-1)} \circ \varphi_{k} \circ \tau_{k}$.
Consider
$\sigma = \lim_{j \to \infty} \tau_{1} \circ \hdots \circ \tau_{j}$
where the limit is considered in the Krull topology. Clearly $\sigma_{j}$
is holomorphic in $W(P) \setminus  E_{j-1}^{\gamma}$
and meromorphic in $W(P)$ for all $j \geq 2$.
Then we define $Y_{P} = (u_{0}^{\gamma} h (\partial{g}/\partial{x}))(0,x_{1},\hdots,x_{n})z \partial / \partial{z}$.
We have $\varphi = {\rm exp}(\sigma_{*} Y_{P})$ since $(a_{1}^{1}z,x_{1},\hdots,x_{n})={\rm exp}(Y_{P})$.
We obtain
\[ \sigma_{*} Y_{P} =
\left({ \left[{ \left({ u_{0}^{\gamma}h \frac{\partial{g}}{\partial{x}} }\right) (0,x_{1},\hdots,x_{n})
\left({ 1 + \sum_{j=2}^{\infty} j \sigma_{j}  {z}^{j-1} }\right) z }\right]
\circ \sigma^{(-1)} }\right)
\frac{\partial}{\partial{z}} . \]
We can change coordinates to obtain
\[ \sigma_{*} Y_{P} =
\left({ u_{0}^{\gamma} +
\sum_{j=1}^{\infty} v_{j}^{P}(x_{1},\hdots,x_{n}) {g}^{j} }\right) g h \frac{\partial}{\partial{x}}  \]
where $v_{j}^{P} \in \vartheta(W(P) \setminus E_{j}^{\gamma})$
and is meromorphic in $W(P)$ for all $j \geq 1$.

 Fix  $Q \in W(P)$ such that
$\pi(Q) \not \in \pi(\overline{\cup_{j \in {\mathbb N}} E_{j}^{\gamma}} \cup (L=1))$.
There exists an open neighborhood $W'$ of $Q$ in $W(P)$ such that
$\pi(\overline{\cup_{j \in {\mathbb N}} E_{j}^{\gamma}} \cup (L=1))$ and $\pi(W')$ are disjoint. This implies that
the series $u_{0}^{\gamma} + \sum_{j=1}^{\infty} u_{j}^{\gamma} {g}^{j}$ and
$u_{0}^{\gamma} + \sum_{j=1}^{\infty} v_{j}^{P} {g}^{j}$ belong to
$\lim_{\leftarrow} G_{W'}/I(W')^{r}$.
By proposition \ref{pro:unqlog}
applied to the fibers of $dx_{1}= \hdots =dx_{n}=0$ we obtain
\begin{equation}
\label{equ:igulog}
 u_{0}^{\gamma} + \sum_{j=1}^{\infty} u_{j}^{\gamma} {g}^{j} =
u_{0}^{\gamma} + \sum_{j=1}^{\infty} v_{j}^{P} {g}^{j} .
\end{equation}
We deduce that $u_{1}^{\gamma} - v_{1}^{P}=0$ contains $W'$. Thus we can extend
the meromorphic function ${(u_{1}^{\gamma})}_{|\gamma}$ to $W(P)$.
By considering every point $P \not \in \gamma \cap (\partial{f} / \partial{x}=0)$ we get
that ${(u_{1}^{\gamma})}_{|\gamma}$ is meromorphic in
\[ (D \cap \gamma) \setminus (E_{1}^{\gamma} \cap (\partial{f} / \partial{x}=0)) = D \cap \gamma  . \]
We can apply the Weierstrass division theorem to obtain that $u_{1}^{\gamma}$ is meromorphic
in $\pi^{-1}(D_{2,\hdots,n+1})$.
By using the equation \ref{equ:igulog} and an
induction process we obtain that ${(u_{j}^{\gamma})}_{|\gamma}$ is meromorphic in
\[ (D \cap \gamma) \setminus (E_{j}^{\gamma} \cap (\partial{f} / \partial{x}=0)) = D \cap \gamma \]
for $j \geq 1$. Thus $u_{j}^{\gamma}$ is meromorphic in $\pi^{-1}(D_{2,\hdots,n+1})$
for all $j \in {\mathbb N}$.
\end{proof}
\begin{pro}
\label{pro:fortra}
  Let $\varphi \in \diff{up}{n+1}$.
Then $\log \varphi$ is t.f. along $Fix \varphi$. Let $\gamma$ be an irreducible component of $Fix \varphi$.
If $\gamma$ is unipotent then $\log \varphi$ is u.t.f. along $\gamma$, otherwise it is
u.s.m. along $\gamma$.
\end{pro}
\begin{proof}
  The results on $\gamma$ are a consequence of propositions \ref{pro:tfuhea} and \ref{pro:pretra}
and lemma \ref{lem:genexp}. Then $\log \varphi$ is t.f. along $Fix \varphi$ by lemma \ref{lem:allforone}.
\end{proof}
We claim that the obstruction for $\log \varphi$ to be u.t.f. along a non-unipotent component $\gamma$
is the existence of germs ${\varphi}_{P}$ for $P \in \gamma$ which can not be embedded in a formal
flow.
\begin{pro}
  Let $\varphi =  {\rm exp}(\hat{u} f \partial/\partial{x}) \in \diff{up}{n+1}$.
Consider a non-unipotent irreducible component $\gamma$ of $Fix \varphi$. Then $\log \varphi$
is u.t.f. along $\gamma$ if and only if there exists a neighborhood of the origin $U$
such that ${\varphi}_{P}$ is embedded in a formal flow for all $P \in \gamma \cap U$.
\end{pro}
\begin{proof}
  The implication $\Rightarrow$ is trivial.
  We denote $L={e}^{u_{0}^{\gamma} h \partial{g}/\partial{x}}$.
We define the set $F= \gamma \cap [(dL \wedge dg=0) \cup (L=1)]$. Since the function $L$ is non-constant in $\gamma$
then $F$ is a proper analytic subset of $\gamma$. We can suppose that the origin belongs to every
irreducible component of $F \cap D$ by shrinking $D$. Moreover we can suppose that
${\varphi}_{Y}$ is embedded in a formal flow for all $Y \in \gamma \cap D$.
Let $P$ be a point in $(D \cap \gamma) \setminus F$. There exists
a system of coordinates $(z,x_{1},\hdots,x_{n})$ in a neighborhood of $P$
(see proof of proposition \ref{pro:pretra}) such that $g=0$ and $\partial{}/\partial{x}$
become $z=0$ and $\partial{} / \partial{z}$ respectively. Consider a neighborhood $W(P)$
of $P$ in $D \cap \gamma$ such that $W(P) \cap F = \emptyset$.

Consider the notations in the proof of proposition \ref{pro:pretra}.
Let $r \in {\mathbb N}$; suppose $a_{j}^{k}$ and $b_{q}$ are holomorphic in $W(P)$
for all $k \leq r$, $j \geq k+1$ and $1 < q \leq r$. We claim that $b_{r+1}$ is
holomorphic in $W(P)$ and then $a_{j}^{r+1} \in \vartheta(W(P))$ for all $j \geq r+2$.
We have $[a_{1}^{1}({(a_{1}^{1})}^{r} - 1)]b_{r+1}=a_{r+1}^{r}$
by the proof of proposition \ref{pro:pretra}. Consider any point $Q \in W(P)$
such that ${(a_{1}^{1}(Q))}^{r} = {L(Q)}^{r} =1$. The diffeomorphism $\varphi_{Q}$
is linearizable by proposition \ref{pro:explin}. There exists an element of
$h \in \diff{}{}$ such that $h^{(-1)} \circ \varphi_{Q} \circ h = a_{1}^{1}(Q)z$
and $j^{r} h = j^{r}((\tau_{1} \circ \hdots \circ \tau_{r-1})_{Q})$ by proposition \ref{pro:concut}.
Now ${(\tau_{1} \circ \hdots \circ \tau_{r-1})}_{Q}^{(-1)} \circ h = z + C_{r+1}(Q) {z}^{r+1} + O({z}^{r+2})$
for some $C_{r+1}(Q) \in {\mathbb C}$. By construction we obtain
\[ [a_{1}^{1}(Q)({(a_{1}^{1}(Q))}^{r} - 1)] C_{r+1}(Q) =a_{r+1}^{r}(Q). \]
Therefore the set $a_{r+1}^{r} = 0$ contains the set ${(a_{1}^{1})}^{r} = 1$.
We have
\[ \emptyset = F \cap W(P) \supset (da_{1}^{1} \wedge dz =0) \cap W(P). \]
Since $a_{1}^{1}$ does not depend on $z$ then $da_{1}^{1}$ never vanishes
in $W(P)$. Therefore any hypersurface $a_{1}^{1}=cte$ is locally irreducible in $W(P)$.
Since $a_{r+1}^{r} \in I(a_{1}^{1} - \mu)$
for every $r$-root of the unit $\mu$ then
\[ b_{r+1} = \frac{a_{r+1}^{r}}{a_{1}^{1}({(a_{1}^{1})}^{r} - 1)}
= \frac{a_{r+1}^{r}}{a_{1}^{1} \prod_{\lambda^{r}=1} (a_{1}^{1} - \lambda)} \]
is holomorphic in $W(P)$. By proceeding as in the proof of proposition
\ref{pro:pretra} we obtain that $u_{j}^{\gamma}$ is holomorphic
in $\pi^{-1}(D_{2,\hdots,n+1} \setminus [\pi(E_{j}^{\gamma}) \cap  \pi(F)])$ for all $j \geq 0$.
Since all the irreducible components of $F$ adhere to $0$ and $0 \not \in \pi(E_{j}^{\gamma})$ then
the codimension of
$\pi(F) \cap \pi(E_{j}^{\gamma}) \subset D_{2,\hdots,n+1}$ is greater or equal than $2$.
By Hartogs' theorem $u_{j}^{\gamma}$ is analytic in $\pi^{-1}(D_{2,\hdots,n+1})$
for all $j \geq 0$. Hence
$\hat{u} = \sum_{j=0}^{\infty} u_{j}^{\gamma} g^{j}$ is u.t.f. along $\gamma$.
\end{proof}
\section{Formal classification}
The goal of this section is providing a complete system of invariants for the formal
classification of up-diffeomorphisms. More precisely, we describe the formal moduli
modulo analytic change of coordinates.
\subsection{Nature of the residue functions}
\label{subsec:narefu}
The formal invariants attached to the germs $\varphi_{P}$ for $P \in Fix \varphi$ are included in
the formal invariants of the up-diffeomorphism $\varphi$. We describe in this subsection the nature
of such invariants.

  Let $\tau \in \diff{u}{}$, there exists a germ of vector field $X \in {\mathcal X}_{N} \cn{}$
such that $X$ is formally conjugated to $\log \tau$. We define the order of contact
$\nu(\tau)$ between $\tau$ and $Id$ as the order of the function $\tau(z)-z$; the definition
does not depend on the choice of the analytical coordinate $z$. We consider the {\it dual form}
of $X$, i.e. the unique form $\omega \in \Omega^{1} \cn{}$ such that $\omega (X)=1$. This form is
meromorphic at the origin.
%
%
We define the residues $Res(\tau)$ and $Res(\log \tau)$ as the residue of $\omega$
at $0$. This residue does not depend on the choice of $X$, moreover it is a formal invariant.
The couple $(\nu(\tau), Res(\tau))$ provides a complete system of formal invariants in $\diff{u}{}$.

 Fix $\varphi \in \diff{up}{n+1}$. Let $\gamma$ be a unipotent non-fibered irreducible component of $Fix \varphi$.
We have $\varphi_{P} \in {\rm Diff}_{u}({\mathbb C},P)$ for all $P \in \gamma$ but
$Q \mapsto Res(\varphi_{Q})$ is not continuous in $\gamma$. An example is provided by
$\varphi={\rm exp}({x}^{2} {(x-y)}^{2} \partial / \partial{x})$, we have that
$Res(\varphi_{(0,y)})=2/{y}^{3}$ if $y \neq 0$ whereas $Res(\varphi_{(0,0)})=0$.
We define
\[ S = \{ P \in \gamma  :  \nu(\varphi_{P}) > \min_{Q \in \gamma} \nu(\varphi_{Q}) \} . \]
The set $S$ is a proper analytic subset of $\gamma$.
The function $Res_{\gamma}(\varphi): \gamma \setminus S \to {\mathbb C}$
defined by $Res_{\gamma}(\varphi)(P) = Res(\varphi_{P})$ is holomorphic. We will prove that
$Res_{\gamma}(\varphi)$ can be extended meromorphically to the whole $\gamma$ as in the example.
Let $\gamma$ be a non-unipotent irreducible component of $Fix \varphi$, we define
$Res_{\gamma}(\varphi): \gamma \to {\mathbb C}$ given by
\[ Res_{\gamma}(\varphi)(Q) = \frac{1}{\ln \left({ \frac{\partial{(x \circ \varphi)}}{\partial{x}} (Q) }\right) } \]
where we choose the determination of $\log$ such that $\log 1 =0$. This definition makes sense since for
$Q = (x_{0},x_{1}^{0}, \hdots, x_{n}^{0}) \not \in \cup_{j=0}^{\infty} E_{j}^{\gamma}$
we have that the residue of
${(\log \varphi)}_{|\cap_{j=1}^{n} (x_{j}=x_{j}^{0})} \in \hat{\mathcal X} ({\mathbb C},x^{0})$
is equal to $Res_{\gamma}(\varphi)(Q)$.
\begin{pro}
\label{pro:restra}
  Let $\varphi = {\rm exp}(\hat{u} f \partial / \partial{x}) \in \diff{up}{n+1}$ and
let $\gamma$ be a non-fibered irreducible component of
$Fix \varphi$. Suppose that $\gamma$ is transversal to $\partial / \partial{x}$ in a neighborhood
of $0$. Then $Res_{\gamma}(\varphi)$ is a meromorphic function of $\gamma$.
\end{pro}
\begin{proof}
Up to a change of coordinates
$\sigma = (x + h(x_{1}, \hdots , x_{n}),x_{1}, \hdots , x_{n})$
we can suppose that $\gamma \equiv (x=0)$. Since $\hat{u}$ is convergent by
restriction to $\gamma$ then $\hat{u}(0,x_{1}, \hdots , x_{n})$
belongs to ${\mathbb C} \{ x_{1} , \hdots , x_{n} \}$.
We can express $f$ in the form $\sum_{j=\nu}^{\infty} f_{j}(x_{1},\hdots,x_{n}) {x}^{j}$
where $\nu = \min_{P \in \gamma} \nu(\varphi_{P})$; in fact we have that
$f_{\nu}(P) \neq 0$ if and only if $\nu(\varphi_{P}) = \nu$ for $P \in \gamma$.
Thus $Res_{\gamma}(\varphi)$ is holomorphic in $\gamma \setminus (f_{\nu}=0)$.

  We consider the transformation
\[ \chi(z,x_{1},\hdots,x_{n})= (z f_{\nu}(x_{1}, \hdots , x_{n}), x_{1},\hdots,x_{n}) . \]
Then $\chi$ is a change of coordinates outside of $f_{\nu}=0$.
We define the diffeomorphism $\tilde{\varphi} = \chi^{(-1)} \circ \varphi \circ \chi$; we obtain
\[ \tilde{\varphi}  =
{\rm exp} \left({ (\hat{u} \circ \chi) {(f_{\nu})}^{\nu}
\left({
z^{\nu} + \sum_{j=\nu+1}^{\infty} f_{j} {(f_{\nu})}^{j-(\nu+1)} {z}^{j}
}\right) \frac{\partial}{\partial{z}} }\right) . \]
We define $\hat{\varphi} = {\rm exp}((\log \tilde{\varphi})/f_{\nu}^{\nu})$.
Since
$(\hat{u} \circ \chi - \hat{u})(0,x_{1},\hdots,x_{n}) \equiv 0$ then the function
$(x_{1},\hdots,x_{n}) \mapsto \nu(\hat{\varphi}_{(0,x_{1},\hdots,x_{n})})$
is equal to the constant $\nu$. Hence the function
$Res_{\gamma}(\hat{\varphi}): \gamma \to {\mathbb C}$ is holomorphic in a neighborhood of
the origin. Note that $f_{\nu}$ does not depend on $z$; thus we obtain
\[ Res_{\gamma}(\varphi)  = Res_{\gamma}(\tilde{\varphi})
=Res_{\gamma}(\hat{\varphi}) / f_{\nu}^{\nu} . \]
Clearly $Res_{\gamma}(\varphi)$ is a meromorphic function of $\gamma$.
\end{proof}
\begin{pro}
\label{pro:resmer}
  Let $\varphi \in \diff{up}{n+1}$ and
let $\gamma$ be a non-fibered irreducible component of
$Fix \varphi$. Then $Res_{\gamma}(\varphi)$ is a meromorphic function of $\gamma$.
\end{pro}
\begin{proof}
Suppose $\gamma$ does not contain orbits of $\partial / \partial{x}$.
As a consequence there exists an irreducible Weierstrass polynomial
\[ H={x}^{k} + a_{k-1}(x_{1},\hdots,x_{n}) {x}^{k-1} + \hdots + a_{0}(x_{1},\hdots,x_{n}) \]
such that $\gamma \equiv (H=0)$. We denote $\pi(x,x_{1},\hdots,x_{n})=(x_{1},\hdots,x_{n})$.
Let $D$ be the critical locus of the projection $\pi_{|\gamma}$.
We denote by $DD$ the critical locus
of the projection $\pi_{|D}$.
The analytic set $DD$ has codimension at least $2$ in $\gamma$.
By proposition \ref{pro:restra} the function $Res_{\gamma}(\varphi)$ is meromorphic in
$\gamma \setminus D$. Suppose that $Res_{\gamma}(\varphi)$ is meromorphic in $\gamma \setminus DD$.
Then there exists a polynomial $R = \sum_{j=0}^{k-1} b_{j}(x_{1},\hdots,x_{n}) {x}^{j}$
such that $Res_{\gamma}(\varphi) \equiv R_{|\gamma}$ by the Weierstrass division.
The coefficients $b_{j}$ ($0 \leq j \leq k-1$) are meromorphic in a neighborhood of the origin
deprived of $\pi(DD)$. Since the codimension of $\pi(DD)$ is at least $2$ in ${\mathbb C}^{n}$ then
$b_{j}$ is meromorphic in a neighborhood of the origin for all $0 \leq j \leq k-1$. Thus
$Res_{\gamma}(\varphi)=R_{|\gamma}$ is meromorphic.

  Let $P \in D \setminus DD$. The set $D$ is transversal to $\partial / \partial {x}$
at $P$ and then smooth. Since $\dim D = n-1$ then there exist coordinates
$(z,z_{1},\hdots,z_{n})$ centered at $P$ such that $\partial / \partial{x}$ and $D$
become $\partial / \partial{z}$ and $(z=0) \cap (z_{1}=0)$ respectively. We consider the
ramification
\[ T(z,z_{1},z_{2}, \hdots ,z_{n})=(z,z_{1}^{k!},z_{2}, \hdots ,z_{n}) . \]
The set $T^{-1}(\gamma)$ is the union of $k$ hypersurfaces, all of them are
transversal to $\partial / \partial{z}$. The function
$Res_{\gamma}(\varphi) \circ T = Res_{T^{-1}(\gamma)}(T^{*} \varphi)$ is meromorphic in
$T^{-1}(\gamma)$ by proposition \ref{pro:restra}.
We undo the ramification to obtain that $Res_{\gamma}(\varphi)$
is meromorphic in a neighborhood of $P$ and then in $\gamma \setminus DD$.

  If $\gamma$ contains orbits of $\partial / \partial{x}$ we can reduce the situation to
the previous one by blowing-up fibered sub-manifolds of $\gamma$. We are done since the field
of meromorphic functions in $\gamma$ is invariant by blow-up.
\end{proof}
Next we express the function $Res_{\gamma}(\varphi)$ in a convenient way.
\begin{lem}
\label{lem:repreteo}
  Let $\varphi \in \diff{up}{n+1}$ and
let $\gamma$ be a non-fibered irreducible component of
$Fix \varphi$. Then there exist $A \in {\mathbb C}\{x,x_{1},\hdots,x_{n}\}$ and
$B \in {\mathbb C}\{x_{1},\hdots,x_{n}\}$ such that $Res_{\gamma}(\varphi) =(A/B)_{|\gamma}$.
\end{lem}
\begin{proof}
There exist $A',B' \in {\mathbb C}\{x,x_{1},\hdots,x_{n}\}$
such that $Res_{\gamma}(\varphi) =(A'/B')_{|\gamma}$ by proposition \ref{pro:resmer}.
We denote $\pi(x,x_{1},\hdots,x_{n})=(x_{1},\hdots,x_{n})$ and $Z= \gamma \cap (B'=0)$.
Let $g=0$ be an irreducible equation of $\gamma$.
We have $\dim Z \leq n-1$ and then $\dim \pi(Z) \leq n-1$. Consider
$h$ in ${\mathbb C}\{x_{1},\hdots,x_{n}\}$ such that $h_{|\pi(Z)} \equiv 0$ but
$h_{|\gamma} \neq 0$; that is possible since $\gamma$ is non-fibered. We obtain
\[ h \in IZ(g,B') \implies h \in \sqrt{(g,B')} \implies h^{m} = Jg + KB' \]
for some $m \in {\mathbb N}$ and $J,K \in  {\mathbb C}\{x,x_{1},\hdots,x_{n}\}$.
The function
\[ \frac{A}{B} = \frac{K A'}{K B' + J g} = \frac{K A'}{{h}^{m}} \]
is equal to $Res_{\gamma}(\varphi)$ in $\gamma$.
\end{proof}
\subsection{The homological equation}
\label{subsec:hoeq}
  We can linearize the problem of formal classification of up-diffeomorphisms.
It can be reduced to deal with equations of the form
\[ \frac{\partial{\alpha}}{\partial{x}} = \frac{A}{f} \]
where $A, f  \in {\mathbb C}\{x,x_{1},\hdots,x_{n}\}$. This equation is called
the {\it homological equation}. Consider the decomposition
$f_{1}^{l_{1}} \hdots f_{p}^{l_{p}} F_{1}^{m_{1}} \hdots F_{q}^{m_{q}}$ of $f$ in irreducible
factors. We suppose that $f_{j}=0$ is non-fibered and $F_{k}=0$ is fibered for
all $1 \leq j \leq p$ and $1 \leq k \leq q$. Moreover we suppose that
$F_{k} \in {\mathbb C}\{x_{1},\hdots,x_{n}\}$ for $1 \leq k \leq q$.
Now we define $f_{N}=\prod_{j=1}^{p} f_{j}^{l_{j}}$ and $f_{F}=\prod_{k=1}^{q} F_{k}^{m_{k}}$.
We say that a meromorphic germ of function $\alpha$ is {\it special} with respect to $f$ if it
can be expressed in the form
\[ \alpha = \frac{\beta}{f_{1}^{l_{1}-1} \hdots f_{p}^{l_{p}-1}  f_{F}} \]
for some $\beta$ in ${\mathbb C}\{x,x_{1},\hdots,x_{n}\}$. A homological equation
$\partial{\alpha}/\partial{x} = A/f$ is {\it special}
if it has a special solution (with respect to $f$). Most of the times we drop the expression
``with respect to f" since it is clear from the context.
We say that the homological equation is
{\it free of residues} if the one-dimensional $1$-form
\[ \left({ \frac{A}{f_{N}}
(x,x_{1}^{0},\hdots,x_{n}^{0}) }\right) dx  \in
{\Omega}_{1}(\cap_{j=1}^{n}(x_{j}=x_{j}^{0})) \]
has residue zero at the neighborhood of every point
$P=(x^{0},x_{1}^{0},\hdots,x_{n}^{0})$ such that
$f_{N}(P)=0$ and $f_{N}(x,x_{1}^{0},\hdots,x_{n}^{0}) \not \equiv 0$.
\begin{lem}
\label{lem:nonevil}
Let $E \equiv [\partial{\alpha} / \partial{x} = A/f]$ be a free of
residues homological equation. Suppose $f_{N}(x,0,\hdots,0) \not \equiv 0$.
Then $E$ is special.
\end{lem}
\begin{proof}
We can suppose that $f_{F}$ is a unit since provided a special solution $\alpha'$ of
$\partial{\alpha'} / \partial{x} = A/f_{N}$ the function
$\alpha = \alpha' /f_{F}$ is a special solution of $E$.
Consider a small neighborhood $V' \times W$ of the origin in ${\mathbb C} \times {\mathbb C}^{n}$
contained in the domain of definition of $E$.
We can suppose that $V'$ is simply connected.
Let $V \subset {\mathbb C}$ a simply connected open set such that $0 \not \in V \subset V'$.
We can suppose that $(V \times W) \cap (f_{N}=0) = \emptyset$.

  The equation $E$ has a solution $\alpha \in \vartheta(V \times W)$.
It can be extended to the set $(V \times W) \setminus (f_{N}=0)$ by
analytic continuation since the residues vanish. Consider a point $Q \in (f_{N}=0)$ such that
$\partial / \partial{x}$ is transversal to $f_{N}=0$ at $Q$. There exist
coordinates $(z,z_{1},\hdots,z_{n})$ centered at $Q$ such that $\partial / \partial{x}$ and $f_{N}=0$
become $\partial / \partial{z}$ and $z=0$ respectively. By integrating with respect
to $z$ we obtain a special solution $\alpha_{Q}$ in the neighborhood of $Q$.
Now $\partial{(\alpha - \alpha_{Q})} / \partial{x}=0$ implies that $\alpha -\alpha_{Q}$
is holomorphic in a neighborhood of $Q$. Hence $\alpha$ is special in a neighborhood of
$Q$. Since $\partial/\partial{x}$ is transversal to $f_{N}=0$ except at a set whose codimension
is greater than $1$ then $\alpha$ is a special solution of $E$.
\end{proof}
\begin{cor}
For $n=1$ (i.e. $A,f  \in {\mathbb C}\{x,x_{1}\}$) we have that a homological equation
$\partial{\alpha} / \partial{x} = A/f$ is free of residues
if and only if it is special.
\end{cor}
In general the vanishing of the residues does not imply the existence of a special solution.
An example is given by $\partial{\alpha} / \partial{x} = 1 /{(z-x y)}^{2}$. On the
one hand it is free of residues since $\alpha = 1/((z-x y)y)$ is a solution.
On the other hand there is no a special solution $\beta /(z-x y)$
since otherwise we obtain  $1 = (\partial{\beta}/\partial{x})(z-x y) + \beta y \in (y,z)$.

  We want to obtain a solution as close as possible to be special for a free of residues equation
$\partial{\alpha} / \partial{x} = A/f$. Let $\prod_{j=1}^{p} f_{j}^{l_{j}}$ be the decomposition
in irreducible factors of $f_{N}$. We consider the analytic set $S(f)$ obtained as
the union of the fibered varieties contained in $\prod_{l_{j}>1} f_{j} = 0$.
The codimension of $S(f)$ is at least $2$. The set $S(f)$ is called {\it evil set}.
\begin{pro}
\label{pro:quaspe}
Let $E \equiv [\partial{\alpha} / \partial{x} = A/f]$ be a free of
residues homological equation. Let $H \in {\mathbb C}\{x_{1},\hdots,x_{n}\} \setminus \{ 0 \}$ vanishing
in the evil set of $f$. Then there exists $k \in {\mathbb N} \cup \{0 \}$ depending only on $f$ such that
$\partial{\alpha} / \partial{x} = (A H^{k})/f$ is special.
\end{pro}
\subsubsection{proof of proposition \ref{pro:quaspe}}
Let $\prod_{j=1}^{p} f_{j}^{l_{j}}$ be the decomposition in irreducible factors of $f_{N}$.
We define the following sheafs:
\begin{itemize}
\item ${\vartheta}_{R}(f)$ of functions $\alpha$ such that
$\alpha f/(f_{1} \hdots f_{p})$ is holomorphic.
\item ${\vartheta}_{D}(f)$ is the sub-sheaf of ${\vartheta}_{R}(f)$ of first integrals of $\partial / \partial{x}$.
\item ${\vartheta}_{Q}(f)$ is the quotient sheaf ${\vartheta}_{R}(f)/{\vartheta}_{D}(f)$.
\end{itemize}
Note that the sheafs ${\vartheta}_{R}(f)$, ${\vartheta}_{D}(f)$ and ${\vartheta}_{Q}(f)$
does not change if we replace $f$ with $f_{F} \prod_{l_{j}>1} f_{j}^{l_{j}}$.

  The vanishing of the residues of $E$ implies that $f_{j}$ divides $A$
if $l_{j}=1$. By replacing $E$ with
$\partial{\alpha}/\partial{x}=(A/\prod_{l_{j}=1} f_{j})/\prod_{l_{j}>1} f_{j}^{l_{j}}$
we can suppose that $l_{j}>1$ for all $1 \leq j \leq p$.

Consider a small polydisk $\Delta \subset {\mathbb C}^{n+1}$
such that $A \in \vartheta(\Delta)$. Now for every
$Y \in \Delta \setminus S(f)$ there exists a solution $\alpha_{Y} \in {\vartheta}_{R}(f)$ of $E$
defined in a neighborhood $U_{Y}$ of $Y$ by lemma \ref{lem:nonevil}.
For another $Y' \in \Delta \setminus S(f)$ we have
$\partial{(\alpha_{Y} - \alpha_{Y'})}/\partial{x} =0$ and then
$\alpha_{Y} - \alpha_{Y'} \in {\vartheta}_{D}(f)(U_{Y} \cap U_{Y'})$.
Therefore $E$ defines a unique section in $H^{0}(\Delta \setminus S(f), \vartheta_{Q}(f))$.
Conversely, let $B \in H^{0}(\Delta \setminus S(f), \vartheta_{Q}(f))$;
we have that $f \partial{B}/\partial{x}$ is holomorphic in $\Delta \setminus S(f)$.
Since ${\rm cod} S(f) \geq 2$ then $f \partial{B}/\partial{x}$ can be extended to $\Delta$.
Thus $\partial{\alpha}/\partial{x} = (f \partial{B}/\partial{x})/f$ is a
free of residues homological equation defined in $\Delta$. We proved that
\begin{lem}
  Fix $f \in {\mathbb C}\{x, x_{1},\hdots,x_{n}\}$.
  The free of residues
homological equations defined in $\Delta$ of the form $\partial{\alpha} / \partial{x} = A/f$
are in a bijective correspondence with $H^{0}(\Delta \setminus S(f), \vartheta_{Q}(f))$.
\end{lem}
The exact sequence
$0 \to \vartheta_{D}(f) \stackrel{i}{\to} \vartheta_{R}(f) \stackrel{p}{\to} \vartheta_{Q}(f) \to 0$
provides the long exact sequence
\[ 0 \to H^{0}(\Delta \setminus S(f), \vartheta_{D}(f)) \stackrel{i^{0}}{\to}
H^{0}(\Delta \setminus S(f), \vartheta_{R}(f)) \stackrel{p^{0}}{\to}  \]
\[ \stackrel{p^{0}}{\to} H^{0}(\Delta \setminus S(f), \vartheta_{Q}(f)) \stackrel{\delta^{0}}{\to}
H^{1}(\Delta \setminus S(f), \vartheta_{D}(f)) \to \hdots . \]
The set of special equations of the form $\partial{\alpha}/\partial{x} = A/f$
is the image of $p^{0}$. We have
\[ 0 \to \frac{H^{0}(\Delta \setminus S(f), {\vartheta}_{Q}(f))}
{p^{0}(H^{0}(\Delta \setminus S(f),{\vartheta}_{R}(f)))}
\stackrel{\delta^{0}}{\rightarrow} H^{1} (\Delta \setminus S(f), {\vartheta}_{D}(f)) . \]

  Since the existence of the evil set is itself an obstruction for the vanishing
of $H^{1} (\Delta \setminus S(f), {\vartheta}_{D}(f))$ we try to remove it by blow-up.
Consider a sequence of blow-ups $\pi_{1}$, $\pi_{2}$, $\hdots$, $\pi_{d}$
centered at fibered varieties. In other words the center of $\pi_{1}$ is
fibered, the center of $\pi_{2}$ is a union of orbits of ${(\pi_{1})}^{*}(\partial / \partial{x})$
and so on. We denote $\pi = \pi_{1} \circ \hdots \circ \pi_{d}$.

  After a finite number of blow-ups we can obtain $\pi$ such that the strict transform $\tilde{f}_{j}$ of
$f_{j}=0$ does not contain orbits of ${\pi}^{*}(\partial / \partial{x})$ for all $1 \leq j \leq p$.
Let $\prod_{j=1}^{q} F_{j}^{m_{j}}$ be the irreducible decomposition of $f_{F}$.
The divisor $[f \circ \pi]$ is of the form
\[ [f \circ \pi] = \sum_{j=1}^{p} l_{j} [\tilde{f}_{j}] + \sum_{j=1}^{q} m_{j} [\tilde{F}_{j}] +
\sum_{j=1}^{s} c_{j} [H_{j}]  \]
where $\tilde{F}_{j}$ is the strict transform of ${F}_{j}=0$. We have that
$\pi^{-1}(S(f)) = \cup_{j=1}^{s} H_{j}$ where
$H_{j}$ is a fibered hypersurface
and $c_{j} \in {\mathbb N}$ for all $1 \leq j \leq s$.
We define $k = \max_{1 \leq j \leq s} c_{j}$. We obtain
\[ \left[{ (H ^{k} f_{F}) \circ \pi }\right] = k [\tilde{H}] +
\sum_{j=1}^{q} m_{j} [\tilde{F}_{j}] + \sum_{j=1}^{s} t_{j} [H_{j}]  \]
where $\tilde{H}$ is the strict transform of $H=0$ and $t_{j} \geq k \geq c_{j}$ for all $1 \leq j \leq s$.

  Consider a polydisk $\Delta$; we consider the sheafs
${\vartheta}_{D}(f \circ \pi)$, ${\vartheta}_{R}(f \circ \pi)$ and ${\vartheta}_{Q}(f \circ \pi)$.
We have
\[ 0 \to \frac{H^{0}(\pi^{-1}(\Delta), {\vartheta}_{Q}(f \circ \pi))}
{p^{0}(H^{0}(\pi^{-1}(\Delta),{\vartheta}_{R}(f \circ \pi)))}
\stackrel{\delta^{0}}{\rightarrow} H^{1} (\pi^{-1}(\Delta) , {\vartheta}_{D}(f \circ \pi)) . \]
The polydisk $\Delta$ is of the form $\Delta_{0} \times \Delta' \subset {\mathbb C} \times {\mathbb C}^{n}$.
Moreover, the choice of $\pi$ implies
$\pi^{-1}(\Delta_{0} \times \Delta') = \Delta_{0} \times \pi^{-1}(\Delta')$. We obtain
\[ H^{1} (\pi^{-1}(\Delta) , {\vartheta}_{D}((H ^{k} f_{F}) \circ \pi)) \sim
H^{1} (\pi^{-1}(\Delta') , \vartheta) . \]
We can prove by using the expression of the blow-up in coordinate charts that
$H^{1} (\pi^{-1}(\Delta') , \vartheta) = H^{1} (\Delta' , \vartheta)$. This implies
\[ H^{1} (\pi^{-1}(\Delta) , {\vartheta}_{D}((H ^{k} f_{F}) \circ \pi)) =0. \]
Since
\[ {\vartheta}_{D}(f \circ \pi) =
{\vartheta}_{D} \left({ \sum_{j=1}^{q} m_{j} [\tilde{F}_{j}] + \sum_{j=1}^{s} c_{j} [H_{j}] }\right)
\subset {\vartheta}_{D} \left({ (H ^{k} f_{F}) \circ \pi }\right)  \]
then for any homological equation $\partial{\alpha} / \partial{x} = A/f$ such that
$A \in {\mathcal O}(\Delta)$ we can find a meromorphic solution $\beta'$ in
${\mathcal O}(\pi^{-1}(\Delta))$ such that
\[ [\beta']_{\infty} \leq
k [\tilde{H}] +  \sum_{j=1}^{p} (l_{j}-1) [\tilde{f}_{j}] +
\sum_{j=1}^{q} m_{j} [\tilde{F}_{j}] + \sum_{j=1}^{s} t_{j} [H_{j}]. \]
By blowing-down we obtain a solution $\beta/(H^{k} f_{F} \prod_{1 \leq j \leq p} {f}_{j}^{l_{j}-1})$
of the equation $\partial{\alpha} / \partial{x} = A/f$ for some $\beta \in {\mathcal O}(\Delta)$.
\subsubsection{Relation between free of residues and special}
Fix $f \in {\mathbb C}\{x,x_{1},\hdots,x_{n}\}$. Let $Fr(f)$ be the set
of free of residues homological equations of the form $\partial{\alpha}/\partial{x}=A/f$
for some $A \in {\mathbb C}\{x,x_{1},\hdots,x_{n}\}$. We define $Sp(f)$ as the subset
of $Fr(f)$ of special equations. We are interested on describing the structure of the
space $Fr(f)/Sp(f)$. We denote
$H^{0}(\Delta \setminus S(f), {\vartheta}_{Q}(f))$ by $Fr(f,\Delta)$ and
$p^{0}(H^{0}(\Delta \setminus S(f),{\vartheta}_{R}(f)))$ by $Sp(f,\Delta)$.
The next proposition is straightforward.
\begin{pro}
Let $0 \neq f \in {\mathbb C}\{x,x_{1},\hdots,x_{n}\}$. We have
\[ \frac{Fr(f)}{Sp(f)} = \lim_{\rightarrow}
\frac{H^{0}(\Delta \setminus S(f), {\vartheta}_{Q}(f))}
{p^{0}(H^{0}(\Delta \setminus S(f),{\vartheta}_{R}(f)))} \]
\end{pro}
We have $Fr(f)/Sp(f)=0$ for small evil sets.
\begin{pro}
Suppose ${\rm cod} S(f) \geq 3$. Then a free of residues $\partial{\alpha}/\partial{x}=A/f$
is special.
\end{pro}
\begin{proof}
It is enough to show that $H^{1} (\Delta \setminus S(f), {\vartheta}_{D}(f)) = 0$ for all
polydisk $\Delta$ small enough. We have
$\vartheta_{D}(f) \sim \vartheta_{D}(f_{N}) = \vartheta$
where $\vartheta$ is the sheaf of holomorphic functions in ${\mathbb C}^{n}$. The sheaf $\vartheta$ is coherent
and then the first homology group does not change by removing
analytic sets of codimension at least $3$ \cite{Scheja}. We obtain
$H^{1} (\Delta \setminus S(f), {\vartheta}_{D}(f)) =0$.
\end{proof}
 Denote by $e(\Delta)$ the canonical mapping from
$Fr(f,\Delta)/Sp(f,\Delta)$ to $Fr(f)/Sp(f)$. The next proposition
makes clear that $Fr(f)/Sp(f)$ behaves like a space of germs.
\begin{pro}
\label{pro:inhoeq}
There exists a fundamental system ${(\Delta_{j})}_{j \in {\mathbb N}}$ of open neighborhoods of the origin
such that $e(\Delta_{j})$ is injective for all $j \in {\mathbb N}$. In particular we have
\[ \frac{Fr(f)}{Sp(f)} = \cup_{j \in {\mathbb N}} \frac{Fr(f, \Delta_{j})}{Sp(f, \Delta_{j})} \]
\end{pro}
\begin{proof}
  We can suppose that $f_{F}$ is a unit without lack of
generality. Therefore we have ${\vartheta}_{D}(f) =\vartheta$ where $\vartheta$ is the sheaf
of holomorphic functions in ${\mathbb C}^{n}$.

  Consider the subset $S'(f)$ of $S(f)$ of points at which $S(f) \subset {\mathbb C}^{n}$
is smooth and of codimension $2$.
We define $\Delta_{j} = V_{j} \times V_{j}^{n}$ such that
\begin{itemize}
\item ${(V_{j})}_{j \in {\mathbb N}}$ is a sequence of open neighborhoods of $0 \in {\mathbb C}$.
\item ${(V_{j}^{n})}_{j \in {\mathbb N}}$ is a sequence of open neighborhoods of $0 \in {\mathbb C}^{n}$.
\item All the connected components of $S'(f)$ in $V_{j}^{n}$ adhere $0$ for $j \in {\mathbb N}$.
\end{itemize}
Fix $j \in {\mathbb N}$. We have to prove that the properties $E \in Fr(f,\Delta_{j})$ and $e(\Delta_{j})(E) =0$
imply $E \in Sp(f,\Delta_{j})$.
Consider the element $\delta^{0}(E)$ of $H^{1} (V_{j}^{n} \setminus S(f), \vartheta)$.
For every open set $U \subset V_{j}^{n} \setminus S(f)$ we define
$\delta^{0}(E,U)$ as the image of $\delta^{0}(E)$ by the canonical mapping
\[ H^{1} (V_{j}^{n} \setminus S(f), \vartheta) \to H^{1} (U  , \vartheta) . \]
We define $S''(f)$ as the subset of $S'(f)$ whose elements $P \in S''(f)$ satisfy
that there exists an  open neighborhood $U_{P} \subset V_{j}^{n}$ of $P$
such that $\delta^{0}(E,U_{P} \setminus S(f))=0$. This property is equivalent
to the existence of a neighborhood of $V_{j} \times \{ P \}$ where $E$ has a special solution.

  By definition $S''(f)$ is open in $S'(f) \cap V_{j}^{n}$. We claim that
$S''(f)$ is closed in $S'(f) \cap V_{j}^{n}$. Let $Q \in \overline{S''(f)} \cap (S'(f) \cap V_{j}^{n})$.
There exists a coordinate system $(y_{1},\hdots,y_{n})$ centered at $Q$ such that
$S(f)$ is given by $y_{1}=y_{2}=0$ in a neighborhood of $Q$.
We denote $K_{\delta}=\cap_{k=1}^{n} (|y_{k}|<\delta)$;
we fix $\delta>0$ such that  $K_{\delta}$ is contained in $V_{j}^{n}$.
Denote $K_{\delta}^{j}=K_{\delta} \setminus (y_{j}=0)$ for $j \in \{1,2\}$; then
$K_{\delta}^{1} \cup K_{\delta}^{2}$
is a Leray covering of $K_{\delta} \setminus S(f)$ for ${\vartheta}$.
There exists a special solution
$h_{k} \in \vartheta_{R}(f)(V_{j} \times K_{\delta}^{k})$ of $E$ for all $k \in \{1,2\}$; moreover
$\delta^{0}(E,K_{\delta} \setminus S(f))$ is given by the function
$h_{1}-h_{2} \in \vartheta(K_{\delta}^{1} \cap K_{\delta}^{2})$.
As a consequence $h_{1}-h_{2}$ can be expressed in the form
\[ h_{1} - h_{2} = \sum_{(k,l) \in {\mathbb Z}^{2}} a_{k,l}(y_{3},\hdots,y_{n}) y_{1}^{k} y_{2}^{l} \]
where $a_{k,l}$ is holomorphic in $\cap_{r=3}^{n} (|y_{r}|<\delta)$ for all $(k,l) \in {\mathbb Z}^{2}$.
Since $Q \in \overline{S''(f)}$ then $a_{k,l}$ vanishes in an open set of $\cap_{r=3}^{n} (|y_{r}|<\delta)$
and then in the whole $\cap_{r=3}^{n} (|y_{r}|<\delta)$ for all $(k,l) \in {({\mathbb Z}_{<0})}^{2}$.
The function
\[ h_{1} - \sum_{k \in {\mathbb Z}, l \geq 0}
a_{k,l}(y_{3},\hdots,y_{n}) y_{1}^{k} y_{2}^{l}=
h_{2} + \sum_{k \geq 0, l<0}
a_{k,l}(y_{3},\hdots,y_{n}) y_{1}^{k} y_{2}^{l} \]
is a special solution of $E$ defined in $V_{j} \times (K_{\delta} \setminus S(f))$
and then in $V_{j} \times K_{\delta}$ since ${\rm cod} S(f) \geq 2$.
Therefore $S''(f)$ is closed in $S'(f) \cap V_{j}^{n}$.

 Every connected component of $S'(f) \cap V_{j}^{n}$ adheres to $0$ and then it intersects $S''(f)$
since $e(\Delta_{j})(E)=0$. We obtain $S''(f)=S'(f) \cap V_{j}^{n}$.
Indeed
$\delta^{0}(E)$ belongs to $H^{1}(V_{j}^{n} \setminus (S(f) \setminus S'(f)), \vartheta)$
by the previous discussion.
The set $S(f) \setminus S'(f)$ has codimension greater than $2$ and then
$\delta^{0}(E)=0$ by Scheja's theorem. Thus $E$ belongs to $Sp(f,\Delta_{j})$.
\end{proof}
In the low dimensional cases we can be even more explicit.
\begin{pro}
Let $0 \neq f \in {\mathbb C}\{x,x_{1},\hdots,x_{n}\}$. Suppose $n \leq 2$.
Then there exists a fundamental system $(\Delta_{j})$ of open neighborhoods of the origin
such that $e(\Delta_{j})$ is an isomorphism. Moreover $Fr(f)/Sp(f)$ is a finite
dimensional complex vector space.
\end{pro}
\begin{proof}
For $n \leq 1$ we have $S(f) = \emptyset$. Thus
$H^{1}(\Delta \setminus S(f), \vartheta_{D}(f))=0$ for every domain $\Delta$ small enough.
That implies $Fr(f)/Sp(f)=0$.

Let $n=2$. We can suppose that $f_{F}$ is a unit
without lack of generality. Moreover we can also suppose that $S(f) \neq \emptyset$
since otherwise we proceed as for $n \leq 1$. Thus we have
$S(f)=\{ (0,0) \}$ since $cod S(f) \geq 2$.
Consider a sequence $\Delta_{j}=B(0,1/j)^{3}$.
We define $K_{j}^{l}=B(0,1/j)^{2} \setminus (x_{l}=0)$ for $l \in \{1,2\}$;
the set $B(0,1/j)^{2} \setminus \{(0,0)\}$ admits a Leray covering $K_{j}^{1} \cup K_{j}^{2}$
for $\vartheta$.

  By proposition \ref{pro:quaspe} there exists $k \in {\mathbb N}$ such that every
$E \in Fr(f,\Delta_{j})$ has a solution $\alpha_{E,l}/x_{l}^{k}$ where
$\alpha_{E,l}$ is special in $B(0,1/j)^{3}$ for $l \in \{1,2\}$.
Now $\delta^{0}(E)$ is given by the function
$\alpha_{E,1}/x_{1}^{k} - \alpha_{E,2}/x_{2}^{k} \in \vartheta(K_{j}^{1} \cap K_{j}^{2})$.
By construction
$\alpha_{E,1}/x_{1}^{k} - \alpha_{E,2}/x_{2}^{k}$ is of the form $h/(x_{1}^{k} x_{2}^{k})$
where $h$ is holomorphic in $B(0,1/j)^{2}$. Since $x_{1}^{a} x_{2}^{b}$
is $0$ in $H^{1}(B(0,1/j)^{2} \setminus \{0\}, \vartheta)$
if $(a,b) \not \in {({\mathbb Z}_{<0})}^{2}$ then the dimension of
$Fr(f,\Delta_{j})/Sp(f,\Delta_{j})$ as a complex vector space is less or equal than $k^{2}$.

The canonical mapping
\[ \frac{Fr(f,\Delta_{j})}{Sp(f,\Delta_{j})} \to \frac{Fr(f,\Delta_{j+1})}{Sp(f,\Delta_{j+1})} \]
is injective by proposition \ref{pro:inhoeq} for $j>>0$. Thus
$\dim_{\mathbb C} Fr(f,\Delta_{j})/Sp(f,\Delta_{j})$ is a non-decreasing sequence
from some moment on. Since it is bounded by above then $\dim_{\mathbb C} Fr(f,\Delta_{j})/Sp(f,\Delta_{j})$
is constant for all $j \geq j_{0}$ and some $j_{0} \in {\mathbb N}$.
Hence $e(\Delta_{j})$ is an isomorphism
for all $j \geq j_{0}$. Clearly $Fr(f)/Sp(f)$ is finite dimensional.
\end{proof}
\subsection{The residue functions are formal invariants}
The invariance of the residues is based on the study of the dual form.
Let $\varphi = {\rm exp}(\hat{u} f \partial/\partial{x}) \in \diff{up}{n+1}$.
We call {\it dual form} of $\log \varphi$ the dual form $dx/(\hat{u}f)$ of
$\log \varphi$ in the relative cohomology of $\partial/\partial{x}$.
Since $\hat{u}$ is t.f. along $f=0$ then there exists $u \in {\mathbb C}\{x,x_{1},\hdots,x_{n}\}$
such that $\hat{u} - u \in (f)$. We define the differential $d_{1} h = (\partial{h}/\partial{x}) dx$
relative to $\partial / \partial{x}$. We have
\[ \frac{dx}{\hat{u} f} = \frac{dx}{uf} + \frac{u-\hat{u}}{f} \frac{1}{u \hat{u}} dx =
 \frac{dx}{uf} + d_{1} \hat{K} \]
for some $\hat{K} \in {\mathbb C}[[x,x_{1},\hdots,x_{n}]]$. Thus the dual form can be decomposed
as the sum of a meromorphic 1-form and a formal exact 1-form. Moreover
$\alpha = {\rm exp}(u f \partial/\partial{x})$ and $\varphi$ have the same residue functions.
Let $f_{N}= \prod_{j=1}^{p} f_{j}^{l_{j}}$;  by lemma \ref{lem:repreteo}
there exist series $P_{j} \in {\mathbb C}\{x,x_{1},\hdots,x_{n}\}$
and $Q_{j} \in {\mathbb C}\{x_{1},\hdots,x_{n}\}$
such that we have $Res_{f_{j}=0}(\varphi)={(P_{j}/Q_{j})}_{|f_{j}=0}$
for all $1 \leq j \leq p$. We define
\[ \omega = \frac{dx}{uf} - \sum_{1 \leq j \leq p} \frac{P_{j}}{Q_{j}} \frac{\partial{f_{j}}/\partial{x}}{f_{j}} dx . \]
The form $\omega$ has vanishing residues in $f_{N}=0$. We obtain
$\omega = d_{1}(A/B)$ for some $A,B \in  {\mathbb C}\{x,x_{1},\hdots,x_{n}\}$
by proposition \ref{pro:quaspe}.  This implies
\begin{lem}
Let $\varphi = {\rm exp}(\hat{u} f \partial/\partial{x}) \in \diff{up}{n+1}$. We have
\[ \frac{dx}{\hat{u} f} = d_{1} \left({ \frac{C}{D} }\right) +
\sum_{1 \leq j \leq p} \frac{P_{j}}{Q_{j}} \frac{\partial{f_{j}}/\partial{x}}{f_{j}} dx  . \]
for some $C \in {\mathbb C}[[x,x_{1},\hdots,x_{n}]]$ and $D \in {\mathbb C}\{x,x_{1},\hdots,x_{n}\}$.
\end{lem}
Let $\varphi_{1}, \varphi_{2} \in \diff{up}{n+1}$. Suppose that there exists
$\hat{\sigma} \in \diffh{}{n+1}$ such that
$\varphi_{2} \circ \hat{\sigma} =\hat{\sigma} \circ \varphi_{1}$.
Denote $f=x \circ \varphi_{1} - x$. Consider the irreducible decomposition
$f_{N}f_{F} = \prod_{j=1}^{p} f_{j}^{l_{j}} \prod_{j=1}^{q} F_{j}^{m_{j}}$ of $f$. We have
\begin{itemize}
\item $\hat{\sigma}$ preserves the fibration $dx_{1} = \hdots = dx_{n}=0$.
\item $f_{j} \circ \hat{\sigma}^{(-1)}=0$ is a non-fibered sub-variety of $Fix \varphi_{2}$
for $1 \leq j \leq p$.
\item $F_{j} \circ \hat{\sigma}^{(-1)}=0$ is a fibered analytic subset of $Fix \varphi_{2}$
for $1 \leq j \leq q$.
\end{itemize}
Let $g_{j}=0$ be an irreducible analytic equation of
$f_{j} \circ \hat{\sigma}^{(-1)}=0$ for $1 \leq j \leq p$. The dual forms of $\log \varphi_{1}$ and
$\log \varphi_{2}$ can be expressed as
\[  d_{1} \left({ \frac{C_{1}}{D_{1}} }\right) +
\sum_{1 \leq j \leq p} \frac{P_{j}^{1}}{Q_{j}^{1}} \frac{\partial{f_{j}}/\partial{x}}{f_{j}} dx
\ \ {\rm and} \ \
 d_{1} \left({ \frac{C_{2}}{D_{2}} }\right) +
\sum_{1 \leq j \leq p} \frac{P_{j}^{2}}{Q_{j}^{2}} \frac{\partial{g_{j}}/\partial{x}}{g_{j}} dx  \]
respectively.
\begin{pro}
\label{pro:rearfor}
The residue functions are formal invariants.
\end{pro}
\begin{proof}
We keep the notations preceding the proposition.
We define $M_{j}=Q_{j}^{1}(Q_{j}^{2} \circ \hat{\sigma})$ and
$N_{j}=P_{j}^{1} (Q_{j}^{2} \circ \hat{\sigma}) - Q_{j}^{1}(P_{j}^{2} \circ \hat{\sigma})$ for $1 \leq j \leq p$.
Our claim is equivalent to $N_{j} \in (f_{j})$
for all $1 \leq j \leq p$.  We denote by $\Omega_{j}$ the dual form of $\log \varphi_{j}$ for $j$ in $\{1,2\}$.
We have $\hat{\sigma}^{*} \Omega_{2} = \Omega_{1}$. Note that
$\hat{\sigma}^{*} (hdx) = (h \circ \hat{\sigma}) (\partial{(x \circ \hat{\sigma})}/\partial{x}) dx$;
we are always working in the relative cohomology with respect to $\partial / \partial{x}$.
We have that $\hat{\sigma}^{*} \Omega_{2} - d_{1}((C_{2}/D_{2})\circ \hat{\sigma})$ is equal to
\[ \sum_{1 \leq j \leq p} \left({ \left({ \frac{P_{j}^{2}}{Q_{j}^{2}}
\frac{\partial{g_{j}}/\partial{x}}{g_{j}} }\right) \circ \hat{\sigma} }\right) d_{1}(x \circ \hat{\sigma}) =
\sum_{1 \leq j \leq p} \frac{P_{j}^{2}}{Q_{j}^{2}} \circ \hat{\sigma}
\frac{\partial{(g_{j} \circ \hat{\sigma})}/\partial{x}}{g_{j} \circ \hat{\sigma}} dx . \]
By construction $g_{j} \circ \hat{\sigma}$ is of the form $\hat{v}_{j} f_{j}$ for some formal
unit $\hat{v}_{j}$ and all $1 \leq  j \leq p$. We obtain
\[ \hat{\sigma}^{*} \Omega_{2} - d_{1}((C_{2}/D_{2})\circ \hat{\sigma}) =
\sum_{1 \leq j \leq p} \frac{P_{j}^{2}}{Q_{j}^{2}} \circ \hat{\sigma}
\frac{\partial{\hat{v}_{j}}/\partial{x}}{\hat{v}_{j}} dx +
\sum_{1 \leq j \leq p} \frac{P_{j}^{2}}{Q_{j}^{2}} \circ \hat{\sigma}
\frac{\partial{f_{j}}/\partial{x}}{f_{j}} dx . \]
The form $\sum_{1 \leq j \leq p} ((P_{j}^{2} / Q_{j}^{2}) \circ \hat{\sigma})
((\partial{\hat{v}_{j}}/\partial{x})/\hat{v}_{j}) dx$ does not have non-fibered poles and
then it can be expressed in the form $d_{1} (J_{1}/J_{2})$ for some
$J_{1} \in {\mathbb C}[[x,x_{1},\hdots,x_{n}]]$ and $J_{2} \in {\mathbb C}[[x_{1},\hdots,x_{n}]]$.
Since $\hat{\sigma}^{*} \Omega_{2}=\Omega_{1}$  then there exist
$C,D \in {\mathbb C}[[x,x_{1},\hdots,x_{n}]]$ such that
\[ d_{1} \left({ \frac{C}{D} }\right) =
\sum_{1 \leq j \leq p} \frac{\partial{f_{j}}/\partial{x}}{f_{j}}
\left({ \frac{P_{j}^{1}}{Q_{j}^{1}} - \frac{P_{j}^{2}}{Q_{j}^{2}} \circ \hat{\sigma} }\right) dx . \]
The previous expression is equivalent to
\[ D^{2} \sum_{j =1}^{p} \frac{\partial{f_{j}}}{\partial{x}} N_{j}
\prod_{k \neq j} f_{k} \prod_{k \neq j} M_{k} =
\prod_{k =1}^{p}  f_{k} \prod_{k=1}^{p} M_{k}
\left({ \frac{\partial{C}}{\partial{x}} D - C \frac{\partial{D}}{\partial{x}}
}\right) . \]
Let $\mu_{r}$ the greatest integer such that $f_{r}^{\mu_{r}}$ divides $D$ for $1 \leq r \leq p$.
For $\mu_{r}>0$ the left hand side belongs to $(f_{r}^{2 \mu_{r}})$ and the right hand side
belongs to $(f_{r}^{\mu_{r}}) \setminus (f_{r}^{\mu_{r}+1})$; this implies that
$\mu_{r}=0$ for $1 \leq r \leq p$. The right hand side belongs to $(f_{r})$, as a consequence
\[ D^{2} \sum_{j =1}^{p} \frac{\partial{f_{j}}}{\partial{x}} N_{j}
\prod_{k \neq j} f_{k} \prod_{k \neq j} M_{k} \in (f_{r}) \implies
N_{r} \in (f_{r}) \]
for all $1 \leq r \leq p$.
\end{proof}
\subsection{Homological equation and formal conjugation}
Let $f$ be an element of ${\mathbb C}\{x,x_{1},\hdots,x_{n}\}$.
We say that $\hat{\sigma} \in \diffh{p}{n+1}$ is {\it special} with respect to $f=0$
if $x \circ \hat{\sigma} - x \in \sqrt{(f_{N})}$.

We say that $\varphi_{1}, \varphi_{2} \in \diff{up}{n+1}$ are {\it formally conjugated by a special
transformation} if there exists a special
$\hat{\sigma} \in \diff{p}{n+1}$ (with respect to $x \circ \varphi_{1}-x=0$)
such that $\varphi_{2} \circ \hat{\sigma} = \hat{\sigma} \circ \varphi_{1}$.
In such a case $(x \circ \varphi_{1} - x)/(x \circ \varphi_{2} -x)$ is a unit.
Thus $\varphi_{1}$, $\varphi_{2}$ both belong to
\[ {\mathcal D}_{f} = \{ \varphi \in \diff{up}{n+1} \ : \ (x \circ \varphi - x)/f \ {\rm is \ a \ unit} \}  \]
for some $f \in {\mathbb C}\{x, x_{1},\hdots,x_{n}\}$ such that $f(0)=0$ and
$(\partial{f}/\partial{x})(0)=0$.
We restrict our study to formal special conjugations. Later on we will see that this point of view is
complete since general formal conjugations can be reduced to the special setting.

 Let $\varphi_{1}, \varphi_{2} \in {\mathcal D}_{f}$. We define
$\hat{u}_{j} =(\log \varphi_{j})(x)/f$ for $j \in \{1,2\}$. Throughout this section we
fix the decomposition $\prod_{j=1}^{p} f_{j}^{l_{j}}$ of $f_{N}$
in irreducible factors. The equation
\[ \frac{\partial{\alpha}}{\partial{x}} = \frac{1}{\hat{u}_{1} f} - \frac{1}{\hat{u}_{2} f} \]
is called the homological equation associated to $\varphi_{1}$ and $\varphi_{2}$. We call it
{\it special} if there exists a {\it special} solution
$\beta/(f_{F} \prod_{j=1}^{p} f_{j}^{l_{j}-1})$
where $\beta$  in  ${\mathbb C}[[x,x_{1},\hdots,x_{n}]]$. Note that
if $1/\hat{u}_{1} - 1/\hat{u}_{2} \in  {\mathbb C}\{x,x_{1},\hdots,x_{n}\}$ then the
definition of special of subsection \ref{subsec:hoeq} implies that the solution is convergent.
Both definitions are the same.
\begin{lem}
\label{lem:undesp}
Consider a homological equation $E \equiv (\partial{\alpha}/\partial{x} = A/f)$ where
$A,f$ belong to ${\mathbb C}\{x,x_{1},\hdots,x_{n}\}$. If there exists a formal special
solution then there also exists a convergent special solution.
\end{lem}
\begin{proof}
We can suppose that $f_{F}$ is a unit without lack of generality.
Consider a formal special solution $\hat{\beta}/(\prod_{1 \leq j \leq p} f_{j}^{l_{j}-1})$ of $E$. We have
\[ \frac{\partial{\hat{\beta}}}{\partial{x}} \prod_{j=1}^{p} f_{j} -
\hat{\beta} \sum_{j=1}^{p} (l_{j}-1) \frac{\partial{f_{j}}}{\partial{x}} \prod_{k \in \{1,\hdots,p\} \setminus \{j\}}
f_{k} = A . \]
If $l_{j}=1$ then $f_{j}$ divides $A$. By considering
$\partial{\alpha}/\partial{x} = (A/\prod_{l_{j}=1} f_{j})/(\prod_{l_{j}>1} f_{j}^{l_{j}})$
we can suppose that $l_{j}>1$ for all $1 \leq j \leq p$. The function $A$ belongs
to the ideal generated by $\prod_{j=1}^{p} f_{j}$ and
$\sum_{j=1}^{p} (l_{j}-1) (\partial{f_{j}} / \partial{x}) \prod_{k \neq j} f_{k}$
in ${\mathbb C}[[x,x_{1},\hdots,x_{n}]]$. Thus it also belongs to the ideal in
${\mathbb C}\{x,x_{1},\hdots,x_{n}\}$ sharing the same generators; in particular
there exist $C,D_{0} \in {\mathbb C}\{x,x_{1},\hdots,x_{n}\}$ such that
\[ A = C \prod_{j=1}^{p} f_{j} - D_{0}
\sum_{j=1}^{p} (l_{j}-1) \frac{\partial{f_{j}}}{\partial{x}} \prod_{k \in \{1,\hdots,p\} \setminus \{j\}} f_{k} . \]
We define $\hat{\gamma} = \hat{\beta} - D_{0}$. We obtain
\[ \frac{\partial{\hat{\gamma}}}{\partial{x}} \prod_{j=1}^{p} f_{j} -
\hat{\gamma} \sum_{j=1}^{p} (l_{j}-1) \frac{\partial{f_{j}}}{\partial{x}} \prod_{k \in \{1,\hdots,p\} \setminus \{j\}}
f_{k} = \left({ C - \frac{\partial{D_{0}}}{\partial{x}} }\right)  \prod_{j=1}^{p} f_{j} . \]
The previous formula implies that $\hat{\gamma} \in (\prod_{j=1}^{p} f_{j})$.
Hence
$(\hat{\gamma}/\prod_{j=1}^{p} f_{j})/(\prod_{j=1}^{p} f_{j}^{l_{j}-2})$ is solution of
\[ \frac{\partial{\alpha}}{\partial{x}} = \frac{C - \partial{D_{0}} / \partial{x}}{\prod_{j=1}^{p} f_{j}^{l_{j}-1}}. \]
By induction on $\max_{1 \leq j \leq p} l_{j}$ we can prove that there exists
$D \in {\mathbb C}\{x,x_{1},\hdots,x_{n}\}$ such that
$(\hat{\beta}-D)/\prod_{j=1}^{p} f_{j}^{l_{j}-1}$ is a solution of an equation
$\partial{\alpha}/\partial{x}=\xi$ for some $\xi \in {\mathbb C}\{x,x_{1},\hdots,x_{n}\}$.
We choose a convergent solution $\alpha_{0} \in {\mathbb C}\{x,x_{1},\hdots,x_{n}\}$ of the latter
equation. Then $D/\prod_{j=1}^{p} f_{j}^{l_{j}-1} + \alpha_{0}$ is a special convergent solution of $E$.
\end{proof}
  We introduce the main proposition in this subsection:
\begin{pro}
\label{pro:hifcon}
Let $\varphi_{1}, \varphi_{2} \in \diff{up}{n+1}$. Then they are formally conjugated by a special
transformation if and only if the associated homological equation is special.
\end{pro}
  The proposition implies that the existence of formal special conjugation is equivalent to
the solvability of a linear differential equation.
The proof of proposition \ref{pro:hifcon} is obtained by reducing the setting to the case
where $\log \varphi_{j}$ is convergent for $j \in \{1,2\}$.
\begin{pro}
\label{pro:cspihsp}
Let ${\varphi}_{j}={\rm exp}(u_{j} f \partial / \partial{x}) \in \diff{up}{n+1}$ with convergent
logarithm for $j \in \{1,2\}$. Assume that the homological equation associated to
$\varphi_{1}$ and $\varphi_{2}$ is special. Then $\varphi_{1}$ and $\varphi_{2}$ are conjugated by a
special diffeomorphism.
\end{pro}
 Denote by $\stackrel{sp}{\sim}$ the equivalence relation given by
$\varphi_{1} \stackrel{sp}{\sim} \varphi_{2}$ if $\varphi_{1}$ is conjugated to $\varphi_{2}$
by a special $\sigma \in \diff{p}{n+1}$.
\begin{proof}
  Let us use the path method (see \cite{Rou:ast} and \cite{Mar:ast}). We define
\[ X_{1+z} = u_{1+z} f \frac{\partial}{\partial{x}} =
\frac{u_{1}u_{2}f}{zu_{1} + (1-z)u_{2}}  \frac{\partial}{\partial{x}} . \]
Denote $c =u_{2}(0)/(u_{2}(0)-u_{1}(0))$. We
have $X_{1+z_{0}} \in {\mathcal X} \cn{n+1}$ for all
$z_{0} \in {\mathbb C} \setminus \{ c \}$. The choice of $X_{1+z}$
assures that the homological equation
\[ \frac{\partial{\alpha}}{\partial{x}} =
z \left({ \frac{dx}{u_{1}f} - \frac{dx}{u_{2}f} }\right) . \]
associated to ${\rm exp}(X_{1})$ and ${\rm exp}(X_{1+z})$ is special.
It is enough to prove that $X_{1} \stackrel{sp}{\sim} X_{2}$
for $c \not \in [0,1]$. If $c \in [0,1]$ we define
\[ Y^{1}_{1+z} = \frac{u_{1}u_{1+i}f}{zu_{1} + (1-z)u_{1+i}}  \frac{\partial}{\partial{x}}
\ \ {\rm and} \ \
Y^{2}_{1+z} = \frac{u_{1+i}u_{2}f}{zu_{1+i} + (1-z)u_{2}}  \frac{\partial}{\partial{x}}. \]
Since $u_{1+i}(0)/(u_{1+i}(0)-u_{1}(0))$ and $u_{1}(0)/(u_{1}(0)-u_{1+i}(0))$ do not belong to
$[0,1]$ then $X_{1} \stackrel{sp}{\sim} X_{1+i} \stackrel{sp}{\sim} X_{2}$.

  Suppose $c \not \in [0,1]$.
We look for $W \in {\mathcal X} \cn{n+2}$ of the form
\[ W = h(x,x_{1},\hdots,x_{n},z) f \frac{\partial}{\partial{x}} + \frac{\partial}{\partial{z}} \]
such that $[W,X_{1+z}]=0$. We ask $h f_{F} \prod_{j=1}^{p} f_{j}^{l_{j}-1}$
to be holomorphic in a connected domain
$V \times V' \subset {\mathbb C}^{n+1} \times {\mathbb C}$ containing $\{ 0 \} \times [0,1]$.
We also require $hf$ to vanish at $\{ 0 \} \times V'$.
Supposed that such a $W$ exists then
${{\rm exp}(W)}_{|z=0}$ is a special mapping conjugating $X_{1}$ and $X_{2}$.

The equation $[W,X_{1+z}]=0$ is equivalent to
\[ u_{1+z} f \frac{\partial{(hf)}}{\partial{x}} - hf \frac{\partial{(u_{1+z}f)}}{\partial{x}}
= \frac{\partial{(u_{1+z} f)}}{\partial{z}} . \]
By simplifying we obtain
\[ u_{1+z} f \frac{\partial{h}}{\partial{x}} - h f \frac{\partial{u_{1+z}}}{\partial{x}} =
\frac{\partial{u_{1+z}}}{\partial{z}}
\Rightarrow \frac{\partial{(h/u_{1+z})}}{\partial{x}} =
 \frac{1}{u_{1} f} - \frac{1}{u_{2} f}. \]
Let $\alpha$ be a special solution of the homological
equation associated to $\varphi_{1}$ and $\varphi_{2}$.
For $p=0$ we can suppose $(\alpha f)(0)=0$ by choosing $\alpha$ of the form
$\alpha'/f_{F}$ where $\alpha' \in {\mathbb C}\{x,x_{1},\hdots,x_{n}\}$
satisfies $\alpha'(0)=0$ and $\partial{\alpha'}/\partial{x}= 1/u_{1}-1/u_{2}$.
We are done by defining $h=u_{1+z} \alpha$.
\end{proof}
The reciprocal is also true.
\begin{pro}
\label{pro:hspicsp}
Let ${\varphi}_{j}={\rm exp}(u_{j} f \partial / \partial{x}) \in \diff{up}{n+1}$ with convergent
logarithm for $j \in \{1,2\}$. Suppose $\varphi_{1} \stackrel{sp}{\sim} \varphi_{2}$.
Then the associated homological equation is special.
\end{pro}
\begin{proof}
  We denote by $T$ the union of the non-fibered irreducible components of $f=0$.
For $P \not \in T$ there exists a special solution $\psi_{j,P}$ of
\[ \frac{\partial \psi_{j,P}}{\partial{x}} = \frac{1}{u_{j} f} \]
defined in a neighborhood of $P$ for $j \in \{1,2\}$. It is a consequence of lemma \ref{lem:nonevil}.
Consider a special diffeomorphism $\sigma$ conjugating $\varphi_{1}$ and $\varphi_{2}$.
We define the diffeomorphism $\sigma_{z} = (1-z)Id + z \sigma$ for all
$z \in {\mathbb C}$; we have that $\sigma_{z}$ is special
for all $z \in {\mathbb C}$. For $P \not \in T$ and $z \in [0,1]$ we define
\[ \gamma_{z}(P) = \psi_{2,P} \circ \sigma_{z} (P) - \psi_{2,P}(P)   . \]
In the previous expression $\psi_{2,P} \circ \sigma_{z} (P)$ is the value at $\sigma_{z}(P)$
of the analytical continuation of $\psi_{2,P}$ along the path $[0,z]:s \to \sigma_{s} (P)$.
Then $\gamma_{z}$ is by construction a special solution of the homological equation
associated to $\sigma_{z}^{(-1)} \circ \varphi_{2} \circ \sigma_{z}$ and $\varphi_{2}$
defined in the complementary of $T$ for all $z \in [0,1]$. There exists a special
solution $\alpha_{P}$ of the homological equation associated to $\varphi_{1}$ and $\varphi_{2}$
and defined in the neighborhood of $P$ for $P \not \in S(f)$ by lemma \ref{lem:nonevil}.
Since $\partial{(\gamma_{1}-\alpha_{P})}/\partial{x}=0$ then $\gamma_{1}$ can be extended to
the complementary of $S(f)$. Moreover ${\rm cod} S(f) \geq 2$ implies that
$\gamma_{1}$ is special in a neighborhood of the origin.
\end{proof}
  Next proposition claims that every formal class of conjugation contains at least one
convergent normal form. That will allow us to prove
proposition \ref{pro:hifcon} by reducing  the problem to the settings
considered in propositions \ref{pro:cspihsp} and \ref{pro:hspicsp}.
\begin{pro}
\label{pro:exfono}
Let $\varphi ={\rm exp}(\hat{u} f \partial / \partial{x}) \in \diff{up}{n+1}$. There exists
a germ of function $u \in {\mathbb C}\{x,x_{1},\hdots,x_{n}\}$ such that $\varphi$ and
${\rm exp}(u f \partial / \partial{x})$ are formally conjugated by a special transformation.
\end{pro}
\begin{proof}
Since $\hat{u}$ is t.f. along $f=0$ (prop. \ref{pro:fortra}) then there exists
$u_{k} \in {\mathbb C}\{x,x_{1},\hdots,x_{n}\}$ such that $\hat{u} - u_{k} \in ({f}^{k})$
for all $k \in {\mathbb N}$.
We denote $\varphi_{k}={\rm exp}(u_{k} f \partial/\partial{x})$.
Let $\gamma_{k}$ be a solution of the homological equation associated to $\varphi_{k}$ and
$\varphi_{1}$. Since $1/u_{k} - 1/u_{1}$ belongs to the ideal $(f)$
we can choose $\gamma_{k}$ in $(x) \cap {\mathbb C}\{x,x_{1},\hdots,x_{n}\}$.
We have $\gamma_{k+1}-\gamma_{k} \in {\mathfrak m}^{k}$ where ${\mathfrak m}$ is the maximal
ideal since $\partial{(\gamma_{k+1} - \gamma_{k})}/\partial{x} \in (f^{k-1})$.
The diffeomorphisms $\varphi_{k}$ and $\varphi_{1}$ are conjugated by
\[ {\sigma}_{k} \stackrel{def}{=} {{\rm exp} \left({ \gamma_{k} \frac{u_{k}u_{1}}{zu_{k} + (1-z)u_{1}}
f \frac{\partial}{\partial{x}} + \frac{\partial}{\partial{z}} }\right)}_{|z=0} . \]
Moreover $\sigma_{k}$ converges in the Krull topology to some special
$\hat{\sigma} \in \diffh{p}{n+1}$ conjugating $\varphi$ and $\varphi_{1}$.
\end{proof}
\begin{proof}[proof of the implication $\Rightarrow$ of proposition \ref{pro:hifcon}]
Define $f=x \circ \varphi_{1} - x$.
We have $\log \varphi_{j} = \hat{u}_{j}f \partial/\partial{x}$ for $j \in \{1,2\}$.
Consider a unit $u_{j}$ in ${\mathbb C}\{x,x_{1},\hdots,x_{n}\}$ such that $\hat{u}_{j}-u_{j} \in (f)$.
Then $\alpha_{j}={\rm exp}(u_{j}f \partial/\partial{x})$ is formally conjugated to $\varphi_{j}$
by a special transformation for $j \in \{1,2\}$ by the proof of proposition \ref{pro:exfono}.
Thus $\alpha_{1}$ and $\alpha_{2}$ are formally conjugated by a special transformation
$\hat{\sigma} \in \diffh{p}{n+1}$. Since
\[ \left({ \frac{1}{\hat{u}_{1}f} - \frac{1}{\hat{u}_{2}f} }\right) -
\left({  \frac{1}{{u}_{1}f} - \frac{1}{{u}_{2}f} }\right) \in {\mathbb C}[[x,x_{1},\hdots,x_{n}]] \]
it is enough to prove that the homological equation associated to $\alpha_{1}$ and $\alpha_{2}$ is special.

 We denote by $\hat{\mathfrak m}$ the maximal ideal of ${\mathbb C}[[x,x_{1},\hdots,x_{n}]]$.
For $k \geq 2$ there exists $h_{k} \in {\mathbb C}\{x,x_{1},\hdots,x_{n}\}$ such that
$h_{k} - x \in (\prod_{j=1}^{p} f_{j})$ and
$(x \circ \hat{\sigma} - h_{k})/ \prod_{j=1}^{p} f_{j} \in {\hat{\mathfrak m}}^{k}$
since $\hat{\sigma}$ is special. We define the special diffeomorphism $\sigma_{k} =(h_{k},x_{1},\hdots,x_{n})$.

Suppose $f_{N}(P) \neq 0$. There exists a special solution
$\psi_{2,P}$ of $\partial{\alpha}/\partial{x}=1/u_{2}f$ defined in the neighborhood of $P$.
We define $\gamma_{k}(P)= \psi_{2,P} \circ \sigma_{k}(P) - \psi_{2,P}(P)$ as in
proposition \ref{pro:hspicsp}. Then $\gamma_{k}$ extends to a special solution defined in
a neighborhood of the origin of the homological equation associated to
$\sigma_{k}^{(-1)} \circ \alpha_{2} \circ \sigma_{k}$ and $\alpha_{2}$ (see proof of prop. \ref{pro:hspicsp}).
We define $\beta_{k} = \gamma_{k} f_{F} \prod_{j=1}^{p} f_{j}^{l_{j}-1}$;
we claim that the sequence $\beta_{k}$ converges to some $\hat{\beta} \in {\mathbb C}[[x,x_{1},\hdots,x_{n}]]$
in the Krull topology. That is a consequence of Taylor's formula since
\[ \gamma_{k} - \gamma_{l}  = \sum_{r=1}^{\infty} \frac{1}{r!}
\left({ \frac{\partial^{r}{\psi_{2}}}{\partial{x^{r}}} \circ \sigma_{l} }\right)
{(x \circ \sigma_{k} - x \circ \sigma_{l})}^{r} \]
implies that $\beta_{k} - \beta_{l} \in \hat{\mathfrak m}^{\min (k,l)}$.
Since $(\log (\sigma_{k}^{(-1)} \circ \alpha_{2} \circ \sigma_{k}))(x)$
converges to $(\log \alpha_{1})(x)$ in the Krull topology
then $\hat{\beta}/(f_{F} \prod_{j=1}^{p} f_{j}^{l_{j}-1})$ is a special
solution of the homological equation associated to $\alpha_{1}$ and $\alpha_{2}$.
\end{proof}

\begin{proof}[proof of the implication  $\Leftarrow$ of proposition \ref{pro:hifcon}]
We have that $\log \varphi_{j}$ is of the form $\hat{u}_{j} f \partial /\partial{x}$
for $j \in \{1,2\}$. Let $u_{j} \in {\mathbb C}\{x,x_{1},\hdots,x_{n}\}$
such that $\hat{u}_{j} - u_{j} \in (f)$ for $j \in \{1,2\}$.
We define $\alpha_{j}= {\rm exp}(u_{j} f \partial / \partial{x})$.
The homological equation
\[ \frac{\partial{\alpha}}{\partial{x}} =
\left({\frac{1}{u_{1}f} - \frac{1}{\hat{u}_{1} f}}\right) +
\left({\frac{1}{\hat{u_{1}} f} - \frac{1}{\hat{u_{2}} f} }\right) +
\left({\frac{1}{\hat{u_{2}} f} - \frac{1}{u_{2}f} }\right)  \]
associated to $\alpha_{1}$ and $\alpha_{2}$ is the result of adding three special equations;
thus it is special. The diffeomorphisms $\alpha_{1}$ and $\alpha_{2}$ are conjugated
by a germ of special diffeomorphism $\sigma$ (proposition \ref{pro:cspihsp}).
By the proof of proposition \ref{pro:exfono} we know that
$\varphi_{j}$ and $\alpha_{j}$ are conjugated by a formal special diffeomorphism $\hat{\sigma}_{j}$
for $j \in \{1,2\}$.
Then $\hat{\sigma}_{2}^{(-1)} \circ \sigma \circ \hat{\sigma}_{1}$ is a formal
special diffeomorphism conjugating $\varphi_{1}$ and $\varphi_{2}$.
\end{proof}
\subsection{Theorem of formal special conjugation}
We describe in this section the nature of the invariants for the formal special conjugation.
The next result is a corollary of  propositions \ref{pro:rearfor} and \ref{pro:hifcon}.
\begin{teo}
The residue functions associated to the non-fibered irreducible components of $f=0$ and
the complex vector space $Fr(f)/Sp(f)$ are a complete system
of formal invariants for the special conjugation in ${\mathcal D}_{f}$.
\end{teo}
Let us clarify the statement. Let $\varphi \in {\mathcal D}_{f}$; we define the
subset ${\mathcal D}_{f}(\varphi)$ of ${\mathcal D}_{f}$
whose elements $\tau$ satisfy that $Res_{\gamma}(\varphi) \equiv Res_{\gamma}(\tau)$ for all
irreducible component $\gamma$ of $f_{N}=0$. The theorem claims the existence of
$Inv_{f}:{\mathcal D}_{f}(\varphi) \to Fr(f)/Sp(f)$ such that
$Inv_{f}(\varphi_{1}) = Inv(\varphi_{2})$ if and only if $\varphi_{1}$ and $\varphi_{2}$
are formally conjugated by a special transformation.

  Consider $\tau={\rm exp}(\hat{u}_{\tau} f \partial / \partial{x}) \in {\mathcal D}_{f}$. Since
$\hat{u}_{\tau}$ is t.f. along $f=0$
there exists a unit $u_{\tau} \in {\mathbb C}\{x,x_{1},\hdots,x_{n}\}$ such that
$\hat{u}_{\tau} - u_{\tau} \in (f)$. We define
a mapping $Inv_{f}^{\varphi}: {\mathcal D}_{f}(\varphi) \to Fr(f)/Sp(f)$ given by
\[ Inv_{f}^{\varphi}(\tau) = \left[{ \frac{\partial{\alpha}}{\partial{x}} =
\frac{1}{f} \left({ \frac{1}{u_{\tau}} - \frac{1}{u_{\varphi}} }\right) }\right] + Sp(f) . \]
The value $Inv_{f}^{\varphi}(\tau)$ is independent of the choices of $u_{\varphi}$ and $u_{\tau}$.
The mapping $Inv_{f}^{\varphi}$ is not the only choice for $Inv_{f}$
since $Inv_{f}^{\varphi} \neq Inv_{f}^{\tau}$ if $Inv_{f}^{\varphi}(\tau) \neq 0$.
Thus $Fr(f)/Sp(f)$ is a classifying space for the formal special conjugation
but the mapping $Inv_{f}$ is not canonical.

  We say that $Fr(f)$ contains units if there exists
$[\partial{\alpha}/\partial{x}=A/f]$ in $Fr(f)$ for some unit
$A \in {\mathbb C}\{x,x_{1},\hdots,x_{n}\}$. Next we prove that there are no redundant
invariants in $Fr(f)/Sp(f)$.
\begin{lem}
The mapping $Inv_{f}^{\varphi}$ is
surjective except if $Fr(f)$ contains units but $Sp(f)$ does not. In such a case
$[\partial{\alpha}/\partial{x}=A/f]+ Sp(f)$ belongs to
$Inv_{f}^{\varphi}({\mathcal D}_{f}(\varphi))$ if and only if $A(0) \neq -1/u_{\varphi}(0)$.
Anyway $Fr(f)/Sp(f)$ is the complex vector space generated by $Inv_{f}^{\varphi}({\mathcal D}_{f}(\varphi))$.
\end{lem}
\begin{proof}
Fix a homological equation $E = [\partial{\alpha}/\partial{x}=A/f] \in Fr(f)$.

Suppose that $1/u_{\varphi}(0) \neq - A(0)$. The formula $1/u = 1/u_{\varphi} + A$
defines a unit
$u \in {\mathbb C}\{x,x_{1},\hdots,x_{n}\}$ such that
$Inv_{f}^{\varphi}({\rm exp}(u f \partial{x})) = E + Sp(f)$.
Then we suppose from now on that $1/u_{\varphi}(0) = - A(0)$.
We have that $\lambda E + Sp(f) \in Inv_{f}^{\varphi}({\mathcal D}_{f}(\varphi))$
for all $\lambda \in {\mathbb C} \setminus \{1\}$.
Note that $\lambda E =  [\partial{\alpha}/\partial{x}=\lambda A/f]$ for $\lambda \in {\mathbb C}$.

Suppose that both $Fr(f)$ and $S(f)$ contain units. Since there exists an equation
$[\partial{\alpha}/\partial{x}=B/f] \in Sp(f)$ such that $B(0) \neq 0$ then
\[ E + Sp(f) = [\partial{\alpha}/\partial{x}=(A+B)/f] + Sp(f) \in Inv_{f}^{\varphi}({\mathcal D}_{f}(\varphi)) . \]

If $Fr(f)$ contains units but $Sp(f)$ does not then
there does not exist a special $[\partial{\alpha}/\partial{x}=B/f]$ such that
$1/u_{\varphi}(0) + A(0) + B(0) \neq 0$. As a consequence
$E+ Sp(f)$ does not belong to $Inv_{f}^{\varphi}({\mathcal D}_{f}(\varphi))$.
\end{proof}
Next results are a direct consequence of the analogous ones on the homological equation
\begin{cor}
Let $f \in {\mathbb C}\{x, x_{1},\hdots,x_{n}\}$. Suppose that $n \leq 1$ or ${\rm cod} S(f) \geq 3$.
Then the residue functions are a complete system of formal invariants for the special conjugation
in ${\mathcal D}_{f} \subset \diff{up}{n+1}$.
\end{cor}
\begin{cor}
Let $f \in {\mathbb C}\{x, x_{1},x_{2}\}$.
Then the residue functions plus a finite number of linear invariants are
a complete system of formal invariants for the special conjugation in
${\mathcal D}_{f} \subset \diff{up}{3}$ .
\end{cor}
\subsubsection{The example $f = {(x_{2}-xx_{1})}^{2} \in {\mathbb C}\{x,x_{1}, \hdots , x_{n} \}$}
\begin{lem}
\label{lem:equide}
Let $f = {(x_{2}-xx_{1})}^{2}$. Consider $E=[\partial{\alpha}/\partial{x}=A/f]$ in $Fr(f)$.
Then $E \in Sp(f)$ if and only if $A \in (x_{1},x_{2})$.
\end{lem}
\begin{proof}
If there exists a special solution $\alpha = \beta/(x_{2}-xx_{1})$ of $E$ then
\[ A=(\partial{\beta}/\partial{x}) (x_{2}-xx_{1}) + \beta x_{1} \in (x_{1},x_{2}). \]
If $A \in (x_{1},x_{2})$ we obtain $A = (x_{2}-xx_{1}) C+ x_{1} D$ for some $C,D$ in
${\mathbb C}\{x,x_{1}, \hdots , x_{n} \}$. Denote
$E'=[\partial{\alpha}/\partial{x} = (C-\partial{D}/\partial{x})/(x_{2}-xx_{1})]$. Since we have
\[ \frac{C-\partial{D}/\partial{x}}{x_{2}-xx_{1}} = \frac{A}{f} - \frac{\partial}{\partial{x}}
\left({ \frac{D}{x_{2}-xx_{1}} }\right) . \]
then $E' \in Fr(f)$. We deduce that
$C - \partial{D}/\partial{x} \in (x_{2}-xx_{1})$; thus $E$ has a special solution
of the form $D/(x_{2}-xx_{1}) + \gamma$ where $\gamma$ is a holomorphic solution of $E'$.
\end{proof}
Consider $E = [\partial{\alpha}/\partial{x} = A/f] \in Fr(f)$.
There exists a solution of $E$ of the form $\beta_{j}/((x_{2}-xx_{1})x_{j})$
for some $\beta_{j} \in {\mathbb C}\{x,x_{1},\hdots,x_{n} \}$ and all $j \in \{1,2\}$
by lemma \ref{lem:equide}.
The element $\delta^{0}(E)$ is given by
$\beta_{1}/((x_{2}-xx_{1})x_{1}) - \beta_{2}/((x_{2}-xx_{1})x_{2})$.
This function is of the form $\beta(x_{1},\hdots,x_{n})/(x_{1}x_{2})$ for some
$\beta \in {\mathbb C}\{x_{1},\hdots,x_{n} \}$.

 Consider $E_{0} = [\partial{\alpha}/\partial{x} = 1/f]$. We have
\[ \frac{\partial}{\partial{x}} \left({ \frac{1}{x_{1}(x_{2}-xx_{1})} }\right) =
\frac{\partial}{\partial{x}} \left({ \frac{x}{x_{2}(x_{2}-xx_{1})} }\right) = \frac{1}{{(x_{2}-xx_{1})}^{2}} . \]
Then $\delta^{0}(E_{0})$ is given by the function $1/(x_{1}x_{2})$. This implies
\begin{lem}
Let $f = {(x_{2}-xx_{1})}^{2} \in {\mathbb C}\{x,x_{1}, \hdots , x_{n} \}$. Then the space
$Fr(f)/Sp(f)$ is equal to ${\mathbb C}\{x_{3},\hdots,x_{n}\} E_{0}$
\end{lem}
Suppose from now on that $n=2$, this implies $Fr(f)/Sp(f) \sim {\mathbb C}$.
Therefore for
$\varphi_{1}= {\rm exp}(\hat{u}_{1} f \partial / \partial{x})$ and
$\varphi_{2}= {\rm exp}(\hat{u}_{2} f \partial / \partial{x})$ such that
$\varphi_{2} \in {\mathcal D}_{f}(\varphi_{1})$ there is a unique
$\lambda \in {\mathbb C}$ such that $1/\hat{u}_{1} - 1/\hat{u}_{2} - \lambda \in (x_{1},x_{2})$.
Then $\varphi_{1}$ and $\varphi_{2}$ are conjugated by a formal special transformation
if and only if $\lambda=0$; this is equivalent to $\hat{u}_{1}(0)=\hat{u}_{2}(0)$.
Note that $2 \hat{u}_{j}(0) = \partial^{2}(x \circ \varphi_{j})/\partial{x_{2}}^{2}(0)$.
We deduce that
\[ \left({ Res_{x_{2}-xx_{1}=0}(\varphi),\frac{\partial^{2}(x \circ \varphi)}{\partial{x_{2}}^{2}}(0,0,0) }\right) \]
is a complete system of formal special invariants in ${\mathcal D}_{f} \subset \diff{up}{3}$.

  We can provide a geometrical interpretation for the non-residual invariant.
Consider $\varphi_{1}={\rm exp}(\hat{u}_{1} f \partial/ \partial{x})$,
$\varphi_{2}={\rm exp}(\hat{u}_{2} f \partial/ \partial{x}) \in {\mathcal D}_{f}$
such that their associated homological
equation $\partial{\alpha}/\partial{x} = A/f$ is free of residues.
Let $P=(x^{0},x_{1}^{0},x_{2}^{0})$ an element of $(f=0) \setminus (x_{1}=x_{2}=0)$.
The one-dimensional germs $\varphi_{1,P}$ and $\varphi_{2,P}$
are formally conjugated by a transformation whose linear part is
$((1-A\hat{u}_{1})(P))w$ where $w$ is a coordinate centered at $P$ in $(x_{1}=x_{1}^{0}) \cap (x_{2}=x_{2}^{0})$.
No other linear part is possible. Since a special conjugation restricted to  $x_{1}=x_{2}=0$
is the identity then the existence of a formal special conjugation at the 1-jet along $f=0$ level
implies $(1 - A \hat{u}_{1}) - 1 \in (x_{1},x_{2}) \Rightarrow A \in (x_{1},x_{2})$.
By lemma \ref{lem:equide} we have that
vanishing of residues plus 1-jet compatibility is equivalent to the existence of
a formal special conjugation.
\section{Convergent actions}
We restricted our study to formal special conjugations. The goal of this section is linking
the equivalence relations ``being formally conjugated" and ``being formally conjugated
by a special transformation". The main result is the following:
\begin{teo}
\label{teo:conact}
  Let $\varphi_{1}, \varphi_{2} \in \diff{up}{n+1}$ be formally conjugated.
Then $\varphi_{1}$ and $\varphi_{2}$ are formally conjugated by a transformation
of the form $\hat{\sigma} \circ \sigma$ where $\sigma$ belongs to $\diff{}{n+1}$ and
$\hat{\sigma} \in \diffh{p}{n+1}$ is special.
\end{teo}
Let us remark that in the theorem $\hat{\sigma}$ is special with respect to
$x \circ \varphi_{2} - x=0$.

  In general a formal conjugation is not of the form $\hat{\sigma} \circ \sigma$.
The action induced by $\hat{\sigma} \circ \sigma$ in the non-fibered
components of $x \circ \varphi_{1} - x=0$ is the one induced by $\sigma$ since $\hat{\sigma}$
is special. Thus such an action is convergent. Now consider
$\varphi_{1}=\varphi_{2}=(x/(1-x),y)$. We have that $\varphi_{1}$ and $\varphi_{2}$
are conjugated by $\hat{\tau}=(x,\sum_{j=1}^{\infty} j! y^{j})$. The action of $\hat{\tau}$
on $x=0$ is not convergent, therefore $\hat{\tau}$ can not be expressed in the form
$\hat{\sigma} \circ \sigma$ where $\sigma \in \diff{}{2}$ and
$\hat{\sigma} \in \diffh{p}{2}$ is special.

In order to prove theorem \ref{teo:conact} it is enough to show
\begin{pro}
\label{pro:exicon}
  Let $\varphi_{j}$ be an element of $\diff{up}{n+1}$ with convergent logarithm
for $j \in \{1,2\}$. Suppose that $\varphi_{1}$ and $\varphi_{2}$ are formally conjugated.
Then they are analytically conjugated.
\end{pro}
Let us explain why proposition \ref{pro:exicon} implies theorem \ref{teo:conact}.
Let $\varphi_{1}, \varphi_{2}$ be elements of $\diff{up}{n+1}$ which are formally conjugated
by $\hat{\tau}$. By proposition \ref{pro:exfono}
there exists $\alpha_{j} \in \diff{up}{n+1}$ such that $\log \alpha_{j}$ is convergent
and $\alpha_{j}$ is conjugated to $\varphi_{j}$ by a special $\hat{H}_{j} \in \diffh{p}{n+1}$
for $j \in \{1,2\}$. We obtain that $\alpha_{1}$ and $\alpha_{2}$
are formally conjugated and then conjugated by some $\tau \in \diff{}{n+1}$.
We define
\[ \hat{\sigma} = \hat{H}_{2} \circ (\tau \circ \hat{H}_{1}^{(-1)} \circ \tau^{(-1)}) \ \ {\rm and} \ \
\sigma = \tau. \]
Clearly $\hat{\sigma} \circ \sigma$ conjugates $\varphi_{1}$ and $\varphi_{2}$. Moreover
$\hat{\sigma}$ is special and $\sigma$ is convergent.

Next proposition is a sort of preparation theorem.
\begin{pro}
\label{pro:prepar}
 Let $\varphi_{1}, \varphi_{2} \in \diff{up}{n+1}$.
Suppose that $\varphi_{1}$ and $\varphi_{2}$ are formally conjugated.
Then for all $\nu \in {\mathbb N}$ there exists $\rho_{\nu}$ in $\diff{}{n+1}$ such that
the diffeomorphism $\varphi_{2,\nu}=\rho_{\nu}^{(-1)} \circ \varphi_{2} \circ \rho_{\nu}$ satisfies
\begin{itemize}
\item $\varphi_{2,\nu} \in {\mathcal D}_{x \circ \varphi_{1} -x} \subset \diff{up}{n+1}$
\item $Res_{\gamma}(\varphi_{1}) \equiv Res_{\gamma}(\varphi_{2,\nu})$
for all non-fibered component $\gamma$ of $Fix \varphi_{1}$.
\item $x \circ \varphi_{1} - x \circ \varphi_{2,\nu} \in {(I(S(x \circ \varphi_{1} - x)) + (x))}^{\nu+1}$.
\end{itemize}
\end{pro}
\begin{proof}
  Denote $f=x \circ \varphi_{1} - x$. Let
$f_{N}f_{F} = \prod_{j=1}^{p} f_{j}^{l_{j}} \prod_{j=1}^{q} F_{j}^{m_{j}}$
be the decomposition of $f$ in irreducible factors. Let $\hat{\rho} \in \diffh{}{n+1}$ be the transformation
conjugating $\varphi_{1}$ and $\varphi_{2}$. Since $\hat{\rho}(Fix \varphi_{1})=Fix \varphi_{2}$
there exist functions
$g_{j} \in {\mathbb C}\{x,x_{1},\hdots,x_{n}\}$ and $G_{k} \in {\mathbb C}\{x_{1},\hdots,x_{n}\}$
such that $(g_{j} \circ \hat{\rho})/f_{j}$ and $(G_{k} \circ \hat{\rho})/F_{j}$ are formal units
for all $1 \leq j \leq p$ and $1 \leq k \leq q$. Consider a function $P_{j}^{1}/Q_{j}^{1}$ such that its
restriction to $f_{j}=0$ is the function $Res_{f_{j}=0}(\varphi_{1})$ for all $1 \leq j \leq p$
(see lemma \ref{lem:repreteo}).
In an analogous way we consider a function
$P_{j}^{2}/Q_{j}^{2}$ such that its
restriction to $g_{j}=0$ is the function $Res_{g_{j}=0}(\varphi_{2})$ for all $1 \leq j \leq p$.
We obtain the system
\[ (X)
\left\{ {
\begin{array}{cc}
g_{j} \circ \hat{\rho} = \hat{v}_{j} f_{j} & {\rm for \ all} \ 1 \leq j \leq p \\
G_{j} \circ \hat{\rho} = \hat{w}_{j} F_{j} & {\rm for \ all} \ 1 \leq j \leq q \\
P_{j}^{1} (Q_{j}^{2} \circ \hat{\rho}) - Q_{j}^{1} (P_{j}^{2} \circ \hat{\rho})
= \hat{r}_{j} f_{j} & {\rm for  \ all} \ 1 \leq j \leq p .
\end{array}
}\right.
\]
The third set of equations is a consequence of the invariance of the residues.
We have that $\hat{v}_{j}$ and $\hat{w}_{k}$ are formal units whereas
$\hat{r}_{j}$ is just a power series.

The ideal $I(S(f))$ associated to the evil set has a
system $L_{1}$, $\hdots$, $L_{d}$ of generators composed by elements of
${\mathbb C}\{x_{1},\hdots,x_{n}\}$. Let us study the equation
$\hat{\rho} \circ \varphi_{1} = \varphi_{2} \circ \hat{\rho}$. The transformation
$\hat{\rho}$ is of the form
\[ \hat{\rho}=\left({
\sum_{j=0}^{2 \nu} \hat{a}_{j} x^{j} + {x}^{2 \nu +1} \hat{A}, \hat{\rho}_{1} , \hdots , \hat{\rho}_{n}
}\right) \]
where $\hat{a}_{0}$, $\hdots$, $\hat{a}_{2 \nu}$, $\hat{\rho}_{1}$, $\hdots$,
$\hat{\rho}_{n} \in {\mathbb C}[[x_{1},\hdots,x_{n}]]$ and $\hat{A} \in {\mathbb C}[[x,x_{1},\hdots,x_{n}]]$.
We denote by $\hat{\rho}'$ be the transformation obtained by replacing $\hat{A}$ with $0$ in the expresion
of $\hat{\rho}$. We want to compare the coefficients of $x^{b}$ ($b \leq \nu$)
of $x \circ \hat{\rho}' \circ \varphi_{1}$ and $x \circ \varphi_{2} \circ \hat{\rho}'$.
We have
\[ \frac{\partial^{b} (x \circ \varphi_{2} \circ \hat{\rho})}{\partial{x}^{b}}(0,x_{1},\hdots,x_{n}) =
\frac{\partial^{b} (x \circ \varphi_{2} \circ \hat{\rho}')}{\partial{x}^{b}}(0,x_{1},\hdots,x_{n}) \]
and
\[ \left({
 \frac{\partial^{b} (x \circ \hat{\rho} \circ \varphi_{1})}{\partial{x}^{b}} -
\frac{\partial^{b} (x \circ \hat{\rho}' \circ \varphi_{1})}{\partial{x}^{b}} }\right) (0,x_{1},\hdots,x_{n})
\in I(S(f))^{\nu+1}
\]
for all $0 \leq b \leq \nu$. The coefficient of $x^{b}$ of $x \circ \hat{\rho}' \circ \varphi_{1}$
can be expressed in the form $C_{b}(x_{1},\hdots,x_{n},\hat{a}_{0}, \hdots , \hat{a}_{2 \nu})$ for some holomorphic
$C_{b}$. Conversely the coefficient of $x^{b}$ of  $x \circ \varphi_{2} \circ \hat{\rho}'$
is of the form $D_{b}(\hat{\rho}_{1} , \hdots , \hat{\rho}_{n},\hat{a}_{0}, \hdots , \hat{a}_{2 \nu})$
for some holomorphic function $D_{b}$. We have
\[ (Y) \ : \  C_{b} - D_{b} = \sum_{k_{1} + \hdots k_{d} = \nu + 1}
\hat{K}_{k_{1},\hdots,k_{d}}  L_{1}^{k_{1}} \hdots L_{d}^{k_{d}}
\ {\rm for \ all} \ 0 \leq b \leq \nu \]
where $\hat{K}_{k_{1},\hdots,k_{d}} \in {\mathbb C}[[x,x_{1},\hdots,x_{n}]]$ for
$k_{1} + \hdots k_{d} = \nu + 1$. By Artin's theorem \cite{Artin}
we can find a solution $(a_{0},\hdots,a_{2 \nu},A,\rho_{1},\hdots,\rho_{n})$
satisfying both the systems $(X)$ and $(Y)$ and such that
\[ \{ a_{0},\hdots,a_{2 \nu},\rho_{1},\hdots,\rho_{n}\} \subset {\mathbb C}\{x_{1},\hdots,x_{n}\}, \
A \in {\mathbb C}\{x,x_{1},\hdots,x_{n}\}. \]
Moreover, we can suppose that $j^{1} a_{k} = j^{1} \hat{a}_{k}$ for all $0 \leq k \leq 2\nu$ and
$j^{1} \rho_{k} = j^{1} \hat{\rho}_{k}$ for all $0 \leq k \leq n$. We define
\[ \rho_{\nu}  = \left({
\sum_{j=0}^{2 \nu} {a}_{j}  x^{j} + x^{2 \nu +1} A, {\rho}_{1}  , \hdots , {\rho}_{n} }\right) . \]
By construction we have that $\rho_{\nu} \in \diff{}{n+1}$ and
$\varphi_{2,\nu} \in {\mathcal D}_{x \circ \varphi_{1} -x}$.
The invariance of the residues by $\rho_{\nu}$ (prop. \ref{pro:rearfor}) implies that
$\varphi_{2,\nu}$ satisfies the second condition in the statement of the proposition.
Now since $(\sum_{j=0}^{2 \nu} {a}_{j}  x^{j} , {\rho}_{1}  , \hdots , {\rho}_{n})$
is a solution of system $(Y)$ then
\[ x \circ \rho_{\nu} \circ \varphi_{1} - x \circ \varphi_{2} \circ \rho_{\nu} \in ({x}^{\nu+1}) + I(S(f))^{\nu+1}. \]
We deduce that $x \circ \varphi_{1} - x \circ \varphi_{2,\nu} \in {((x) + I(S(f)))}^{\nu+1}$.
\end{proof}
We intend to prove that the homological equation associated to $\varphi_{1}$ and $\varphi_{2,\nu}$
is special for $\nu >>0$. We will construct special solutions in the neighborhood of every point
outside of a set of codimension greater or equal than $3$.
Next lemmas are of technical interest.
\begin{lem}
\label{lem:arre1}
Let $I,J$ be ideals of a noetherian ring $R$.
Then there exists $\nu_{0} \in {\mathbb N}$ such that we have
$J^{\nu} \cap I \subset J^{\nu-\nu_{0}}I$ for all $\nu \geq \nu_{0}$.
\end{lem}
\begin{proof}
The equation $J^{\nu} \cap I = J^{\nu-\nu_{0}}(J^{\nu_{0}} \cap  I)$ is a consequence of
Artin-Rees lemma (see corollary 10.10 in \cite{At-Mc}); this implies $J^{\nu} \cap I \subset J^{\nu-\nu_{0}}I$.
\end{proof}
\begin{lem}
\label{lem:arre2}
Let $R$ be a domain of integrity. Consider an element $g$ in $R \setminus \{ 0 \}$ and an
ideal $J \subset A$. Then there exists $\nu_{0} \in {\mathbb N}$ such that
$(J^{\nu}:g) \subset J^{\nu-\nu_{0}}$ for all $\nu \geq \nu_{0}$.
\end{lem}
\begin{proof}
We define $I=(g)$. Consider $h \in (J^{\nu}:g)$, we have $hg \in I \cap J^{\nu}$.
By lemma \ref{lem:arre1} there exists $\nu_{0} \in {\mathbb N}$ such that
$J^{\nu} \cap I \subset J^{\nu-\nu_{0}}I$ for $\nu \geq \nu_{0}$. Therefore
$h$ belongs to $J^{\nu-\nu_{0}}$.
\end{proof}
Let $I$ be an ideal of ${\mathbb C}\{x,x_{1},\hdots,x_{n}\}$.
Fix a finite system $L=\{ L_{1},\hdots,L_{d} \}$ of generators of $I$.
There exists a neighborhood $W_{L}$ of the origin such that
$L_{j} \in \vartheta(W_{L})$ for $1 \leq j \leq d$. For $P=(x^{0},x_{1}^{0}, \hdots , x_{n}^{0}) \in W_{L}$ we
define the ideal $I_{P}$ contained in the ring
${\mathbb C}\{x-x^{0},x_{1}-x_{1}^{0},\hdots,x_{n}-x_{n}^{0}\}$ and generated by
$L_{1},\hdots,L_{d}$. The definition of $I_{P}$ does not depend on the
system of generators.
For a different finite system of generators $L'$ there exists
a neighborhood of the origin $W_{LL'} \subset W_{L} \cap W_{L'}$ where both definitions of
$I_{P}$ coincide for all $P \in W_{LL'}$.
\begin{lem}
\label{lem:germ}
 Let $0 \neq f \in {\mathbb C}\{x,x_{1},\hdots,x_{n}\}$.
There exists $\nu_{0} \in {\mathbb N}$ such that for all $\nu \geq \nu_{0}$
we have an open set $U_{\nu} \ni 0$ satisfying that
for all $P \in U_{\nu}$ and $A$ in $((x)+I(S(f)))_{P}^{\nu}$
such that $\partial{\alpha}/\partial{x}=A/f$ is special in a neighborhood of $P$
there exists a special solution
$\beta_{P}/(f_{F} \prod_{j=1}^{p} f_{j}^{l_{j}-1})$ with
$\beta_{P} \in ((x)+I(S(f)))_{P}^{\nu-\nu_{0}}$.
\end{lem}
\begin{proof}
Let $f_{N} = \prod_{j=1}^{p} f_{j}^{l_{j}}$.
Denote $J=(x) + I(S(f))$.
The proof is by induction on $l=\max_{j=1}^{p} l_{j}$.
If $l =0$ we can choose $\beta \in J^{\nu+1}$.

There exists $\nu_{1} \in {\mathbb N}$ such that $(J^{\nu}:f_{j}) \subset J^{\nu-\nu_{1}}$
for all $\nu \geq \nu_{1}$ and all $1 \leq j \leq p$ by lemma \ref{lem:arre2}. Indeed we obtain
$(J_{P}^{\nu}:f_{j}) \subset J_{P}^{\nu-\nu_{1}}$ for all $P$ in some open set $U_{\nu}^{1} \ni 0$
by Oka's coherence theorem (see \cite{Gunning} page 67).

 Denote $E=[\partial{\alpha}/\partial{x}=A/f]$ and
$U_{\nu}^{2} = \cap_{j=1}^{p} U_{\nu - (j-1)\nu_{1}}^{1}$.
We can suppose $l_{j} \neq 1$ for $1 \leq j \leq p$ since
otherwise we replace $E$ with
\[ \frac{\partial{\alpha}}{\partial{x}} =
\frac{A/\prod_{l_{j}=1} f_{j}}{f_{F} \prod_{l_{j} \neq 1} f_{j}^{l_{j}}} \]
where $A/\prod_{l_{j}=1} f_{j} \in J_{P}^{\nu - p \nu_{1}}$ for all
$P \in  U_{\nu}^{2}$ and all $\nu \geq p \nu_{1}$. A special solution
$\beta'/(f_{F} \prod_{j=1}^{p} f_{j}^{l_{j}-1})$ of $E$ is characterized by
\begin{equation}
\label{equ:ide1}
\frac{\partial{\beta'}}{\partial{x}} \prod_{j=1}^{p} f_{j} -
\beta' \sum_{j=1}^{p} (l_{j}-1) \frac{\partial{f_{j}}}{\partial{x}} \prod_{k \in \{1,\hdots,p\} \setminus \{j\}}
f_{k} = A .
\end{equation}
We define the ideal
\[ I=
\left({ \prod_{j=1}^{p} f_{j}, \sum_{j=1}^{p} (l_{j}-1) \frac{\partial{f_{j}}}{\partial{x}}
\prod_{k \in \{1,\hdots,p\} \setminus \{j\}} f_{k} }\right) . \]
Since $E$ is special in the neighborhood of $P$ then $A \in I_{P}$.
By lemma \ref{lem:arre1} and Oka's theorem there exist $\nu_{3} \in {\mathbb N}$
such that $I_{P} \cap J_{P}^{\nu} \subset I_{P} J_{P}^{\nu - \nu_{3}}$ for all $P$ in
some open set $U_{\nu}^{3} \ni 0$ and all $\nu \geq \nu_{3}$. As a consequence there exist
$B_{P} \in J_{P}^{\nu-p\nu_{1}-\nu_{3}}$ and $C_{P} \in J_{P}^{\nu-p\nu_{1}-\nu_{3}}$
such that
\begin{equation}
\label{equ:ide2}
B_{P}  \prod_{j=1}^{p} f_{j} -
C_{P} \sum_{j=1}^{p} (l_{j}-1) \frac{\partial{f_{j}}}{\partial{x}} \prod_{k \in \{1,\hdots,p\} \setminus \{j\}}
f_{k} = A .
\end{equation}
for all $P \in U_{\nu}^{2} \cap U_{\nu-p\nu_{1}}^{3}$ and all $\nu \geq p \nu_{1} + \nu_{3}$.
By subtracting the equations \ref{equ:ide1} and \ref{equ:ide2} we obtain
$\beta' - C_{P} \in (\prod_{j=1}^{p} f_{j})$. Therefore the function
\[ \eta = \frac{(\beta'-C_{P})/\prod_{j=1}^{p} f_{j}}{f_{F} \prod_{j=1}^{p} f_{j}^{l_{j}-2}} \]
is a special solution of
\begin{equation}
\label{equ:reheq}
 \frac{\partial{\alpha}}{\partial{x}} =
\frac{B_{P} - \partial{C_{P}} / \partial{x}}{f_{F} \prod_{j=1}^{p} f_{j}^{l_{j}-1}}
\end{equation}
(see proof of lemma \ref{lem:undesp}). We have that
$B_{P} - \partial{C_{P}}/\partial{x} \in J_{P}^{\nu-p\nu_{1}-\nu_{3}-1}$
for all $\nu \geq p\nu_{1} + \nu_{3} + 1$. By the hypothesis of induction there exists
$\nu_{4} \in {\mathbb N}$ and a special solution
$\gamma_{P}/(f_{F} \prod_{j=1}^{p} f_{j}^{l_{j}-2})$
of equation \ref{equ:reheq} such that $\gamma_{P}$ belongs to $J_{P}^{\nu-p\nu_{1}-\nu_{3}-1-\nu_{4}}$
for all $\nu \geq p\nu_{1}+\nu_{3}+\nu_{4}+1$ and all $P \in U_{\nu}^{4}$ for some open set
$U_{\nu}^{4} \ni 0$. We define $\nu_{0} = p\nu_{1}+\nu_{3}+\nu_{4}+1$
and $U_{\nu} = U_{\nu}^{2} \cap U_{\nu-p\nu_{1}}^{3} \cap U_{\nu}^{4}$.
The function $\beta_{P} = C_{P} + \gamma_{P} \prod_{j=1}^{p} f_{j}$ belongs to $J_{P}^{\nu-\nu_{0}}$.
Moreover $\beta_{P}/(f_{F} \prod_{j=1}^{p} f_{j}^{l_{j}-1})$ is a special
solution of $E$ in the neighborhood of $P$ if $P \in U_{\nu}$.
\end{proof}
  Next we prove proposition \ref{pro:exicon}. The proof is based on the fact that
in the neighborhood of the generic points of $S(f)$ the quotient
``free of residues homological equations / special equations" generates a finite dimensional
vector space over the meromorphic functions in $S(f)$.
\begin{proof}[proof of proposition \ref{pro:exicon}]
Let $\varphi_{2,\nu}$ be the diffeomorphism and $U_{\nu}$ be the open set
given by proposition \ref{pro:prepar} for
all $\nu \in {\mathbb N}$. It is enough to prove that there exists $\nu_{0} \in {\mathbb N}$
such that $\varphi_{1} \stackrel{sp}{\sim} \varphi_{2, \nu}$ for all $\nu \geq \nu_{0}$.
Denote $f= x \circ \varphi_{1} -x$.
Let $f_{N} = \prod_{j=1}^{p} f_{j}^{l_{j}}$.
We have that $\varphi_{1}$ and $\varphi_{2,\nu}$ are
of the form ${\rm exp}(u_{1} f \partial / \partial{x})$ and ${\rm exp}(u_{2,\nu} f \partial / \partial{x})$
respectively. Consider the homological equation $E_{\nu} = [\partial{\alpha}/\partial{x} = A_{\nu}/f]$
associated to $\varphi_{1}$ and $\varphi_{2,\nu}$. The equation $E_{\nu}$ is free of residues
by proposition \ref{pro:prepar}. Denote $J=(x)+I(S(f))$;
we claim that $u_{1}f -u_{2,\nu}f \in J^{\nu+1}$.
Otherwise we have $u_{1}f - u_{2,\nu}f \in J^{a} \setminus J^{a+1}$ for some
$a < \nu+1$. Note that since $S(f) \subset (\prod_{l_{j} \geq 2} f_{j}=0)$ then
$f \in I(S(f))^{2}$; this implies that $f \in J^{2}$ and $a \geq 2$. This property can be used to prove
that
\[ (u_{1} f \partial/\partial{x})^{j} (x) - (u_{2,\nu} f \partial/\partial{x})^{j}(x) \in J^{a+1} \]
for all $j \geq 2$. As a consequence we obtain that
$x \circ \varphi_{1} - x \circ \varphi_{2,\nu} \not \in J^{a+1}$ and that is impossible since $a+1 \leq \nu+1$.
Since $A_{\nu}=(u_{2,\nu}f - u_{1}f)/(u_{1} u_{2,\nu} f)$ then $A_{\nu} \in J^{\nu-\nu_{1}}$
for all $\nu \geq \nu_{1}$ and some $\nu_{1} \in {\mathbb N}$. The function $A_{\nu}$ is defined in
some open set $U_{\nu}^{1} \ni 0$.

  Let $T(f)$ be the set of points of $S(f)$ where $S(f)$ is smooth and of local codimension $2$.
Consider $P=(0,x_{1}^{0},\hdots,x_{n}^{0}) \in T(f)$; there exists $k(P) \in {\mathbb N}$ such that
\[ \frac{\partial{\alpha}}{\partial{x}} =  \frac{H^{k(P)} A_{\nu}}{f} \]
is special in the neighborhood of $P$
for every $H \in {\mathbb C}\{x_{1}-x_{1}^{0},\hdots,x_{n}-x_{n}^{0}\}$ vanishing in $S(f)$
and all $\nu \in {\mathbb N}$. Moreover, a review of the proof of proposition \ref{pro:quaspe}
implies that we can choose the same $k=k(P)$ for all $P \in (x=0) \cap (T(f) \setminus F(f))$
where $F(f) \subset S(f)$ is a fibered analytic variety such that ${\rm cod} (F(f)) \geq 3$.

Fix $P= (0,x_{1}^{0},\hdots,x_{n}^{0}) \in T(f) \setminus F(f)$. We can find
new coordinates $(y_{1},\hdots,y_{n})$ centered at
$(x_{1},\hdots,x_{n})=(x_{1}^{0},\hdots,x_{n}^{0})$ such that $S(f)=(y_{1}=0) \cap (y_{2}=0)$.
Suppose that $P \in U_{\nu-\nu_{1}} \cap U_{\nu}^{1}$;
by lemma \ref{lem:germ} there exists $\nu_{2} \in {\mathbb N}$ such that the equation $E_{\nu}$ has a solution
\[ \alpha_{\nu,P,j} =
\frac{\beta_{\nu,P,j}}{y_{j}^{k} f_{F} \prod_{r=1}^{p} f_{r}^{l_{r}-1}} \]
where $\beta_{\nu,P,j} \in J_{P}^{\nu-\nu_{1}-\nu_{2}}$ for all $\nu \geq \nu_{1}+\nu_{2}$,
and $j \in \{1,2\}$.
Consider the set
$K_{\delta}=(|x|<\delta) \cap \cap_{j=1}^{n} (|y_{j}| < \delta)$ for some $\delta=\delta(\nu, P)>0$
small enough.
The element $\delta^{0}(E_{\nu},K_{\delta} \setminus S(f))$ of
$H^{1}(K_{\delta} \setminus S(f), \vartheta_{D}(f))$
is given by the function $\alpha_{\nu,P,1} - \alpha_{\nu,P,2}$. We obtain
\[  (\alpha_{\nu,P,1} - \alpha_{\nu,P,2}) f_{F} =
\frac{y_{2}^{k} \beta_{\nu,P,1}  - y_{1}^{k} \beta_{\nu,P,2}}{y_{1}^{k} y_{2}^{k} \prod_{j=1}^{p} f_{j}^{l_{j}-1}} . \]
Since $\partial{(\alpha_{\nu,P,1} - \alpha_{\nu,P,2})}/\partial{x}=0$
then $\prod_{j=1}^{p} f_{j}^{l_{j}-1}$ divides $y_{2}^{k} \beta_{\nu,P,1}  - y_{1}^{k} \beta_{\nu,P,2}$.
Hence there exists an open set $U_{\nu}^{2} \ni 0$ such that for $P \in U_{\nu}^{2}$ the function
$(\alpha_{\nu,P,1} - \alpha_{\nu,P,2})f_{F}$ can be expressed in the form
$h_{\nu,P}(y_{1},\hdots,y_{n})/(y_{1}^{k}y_{2}^{k})$ where $h_{\nu,P} \in J_{P}^{\nu-\nu_{1}-\nu_{2}-\nu_{3}}$
for all $\nu \geq \nu_{1} + \nu_{2} + \nu_{3}$ and some $\nu_{3} \in {\mathbb N}$.
We define $\nu_{0} = \nu_{1} + \nu_{2} + \nu_{3} + (2k-1)$.
The set $J_{P}^{b} \cap {\mathbb C} \{y_{1},\hdots,y_{n}\}$ is contained in
$I(S(f))_{P}^{b} \subset (y_{1},y_{2})^{b}$ for all $b \in {\mathbb N}$,
and then for $\nu \geq \nu_{0}$ the function $h_{\nu,P} /(y_{1}^{k}y_{2}^{k})$ is of the form
\[ \frac{h_{\nu,P}}{y_{1}^{k}y_{2}^{k}} = H + \sum_{0 < j \leq k} \frac{c_{j}(y_{2},\hdots,y_{n})}{ y_{1}^{j}} +
\sum_{0 < j \leq k} \frac{d_{j}(y_{1},y_{3},\hdots,y_{n})}{ y_{2}^{j}} \]
where $H$, $c_{j}$ and $d_{j}$ are holomorphic in $K_{\delta}$ for all $0 <j \leq k$. The function
\[
\alpha_{\nu,P}  \stackrel{def}{=}
\alpha_{\nu,P,1} - (H + \sum_{-k \leq j < 0} c_{j}  y_{1}^{j})/f_{F}
= \alpha_{\nu,P,2} + (\sum_{-k \leq j < 0} d_{j} y_{2}^{j})/f_{F} \]
is a special solution of $E_{\nu}$ in $K_{\delta}$ for $\nu \geq \nu_{0}$.

Consider a polydisk
$0 \in \Delta_{\nu}$ in the variables $(x,x_{1},\hdots,x_{n})$ contained in the set
$U_{\nu-\nu_{1}} \cap U_{\nu}^{1} \cap U_{\nu}^{2}$.
We denote $\pi(x,x_{1},\hdots,x_{n})=(x_{1},\hdots,x_{n})$. For all $\nu \geq \nu_{0}$ and
$P \in (x=0) \cap (T(f) \setminus F(f)) \cap \Delta_{\nu}$
there exists a polydisk $\Delta_{\nu,P} \subset \Delta_{\nu}$ centered at $P$
and a special solution $\alpha_{\nu,P}$ of $E_{\nu}$ defined in $\Delta_{\nu,P}$. By using the homological
equation we can extend $\alpha_{\nu,P}$ to $\Delta_{\nu} \cap \pi^{-1}(\Delta_{\nu,P} \setminus S(f))$
and then to $\Delta_{\nu} \cap \pi^{-1}(\Delta_{\nu,P})$ since ${\rm cod} S(f) \geq 2$. We obtain
special solutions of $E_{\nu}$ ($\nu \geq \nu_{0}$) in the neighborhood of every
point not in $(S(f) \setminus T(f)) \cup F(f)$. Therefore we have
\[ \delta^{0}(E_{\nu}) \in H^{1}(\Delta_{\nu} \setminus [(S(f) \setminus T(f)) \cup F(f)],\vartheta_{D}(f)) . \]
Since the codimension of $(S(f) \setminus T(f)) \cup F(f)$ is greater or equal than $3$ then
$\delta^{0}(E_{\nu})=0$ for $\nu \geq \nu_{0}$. We deduce that $E_{\nu} \in Sp(f)$ and then
$\varphi_{1} \stackrel{sp}{\sim} \varphi_{2,\nu}$ for $\nu \geq \nu_{0}$ by proposition \ref{pro:cspihsp}.
\end{proof}
\bibliography{rendu}
\end{document}